%% file: main.tex
\documentclass{article}
\usepackage{hyperref}
\usepackage{amssymb}
\usepackage{amsthm}
\usepackage{microtype}
\usepackage{array}
\usepackage{xspace}
\usepackage{csquotes}
\usepackage{mathtools}
\usepackage{tikz-cd}
\usepackage{bbm}
\usepackage[numbers,sort]{natbib}
\usepackage{enumerate}
\usepackage{todonotes}
\usepackage{enumitem}

\newtheorem{thm}{Theorem}

\newtheorem{lemma}[thm]{Lemma}

\newtheorem{corollary}[thm]{Corollary}
\newtheorem{proposition}[thm]{Proposition}
\theoremstyle{definition}
\newtheorem{defi}[thm]{Definition}
\newtheorem{definition}[thm]{Definition}
\theoremstyle{remark}
\newtheorem{example}[thm]{Example}
\newtheorem{remark}[thm]{Remark}

\newcommand{\N}{\ensuremath{\mathbb{N}}}
\newcommand{\Z}{\ensuremath{\mathbb{Z}}}

\newcommand{\such}{\ensuremath{\;|\;}} 
\DeclarePairedDelimiter{\set}{\{}{\}}

\DeclareMathOperator{\swap}{swap}
\DeclareMathOperator{\Set}{Set}
\DeclareMathOperator{\Cat}{Cat}
\DeclareMathOperator{\Aut}{Aut}

\DeclareMathOperator{\op}{op}

\newcommand{\wcat}{\mathbb}
\newcommand{\coslice}[2]{{\scriptstyle #2}\backslash #1}
\DeclareMathOperator{\id}{id}

\DeclareMathOperator{\inc}{in}

\DeclareMathOperator{\pr}{pr}
\DeclareMathOperator{\ind}{ind}
\newcommand{\G}{\ensuremath{\mathbb{G}}}
\DeclareMathOperator{\List}{List}
\DeclareMathOperator{\bat}{Bat}
\DeclareMathOperator{\glob}{Glob}
\newcommand{\bipointed}{\glob^{\star\star}}
\newcommand{\compbip}{\comp^{\star\star}}
\newcommand{\catbip}{\omega\Cat^{\star\star}}
\newcommand{\freebip}{\free^{\star\star}}
\newcommand{\cellbip}{\cell^{\star\star}}
\newcommand{\Tbip}{T^{\star\star}}
\DeclareMathOperator{\src}{src}
\DeclareMathOperator{\tgt}{tgt}
\DeclareMathOperator{\btlist}{br}
\newcommand{\bt}[1]{\btlist[#1]}
\DeclareMathOperator{\susp}{\Sigma}
\DeclareMathOperator{\deloop}{\Omega}
\DeclareMathOperator{\PosOp}{Pos}
\newcommand{\Pos}[2][]{\PosOp_{#1}(#2)}

\newcommand{\disk}[1]{\ensuremath{\mathbb{D}^{#1}}}

\DeclareMathOperator{\comp}{Comp}
\DeclareMathOperator{\cell}{Cell}
\DeclareMathOperator{\type}{Sph}
\DeclareMathOperator{\ty}{bdry}
\DeclareMathOperator{\Par}{Par}
\DeclareMathOperator{\free}{Free}

\DeclareMathOperator{\supp}{supp}
\DeclareMathOperator{\str}{str}
\newcommand{\bdry}[2][]{\ensuremath{\partial_{#1}#2}}
\newcommand{\srcps}[2][]{\ensuremath{s_{#1}^{#2}}}
\newcommand{\tgtps}[2][]{\ensuremath{t_{#1}^{#2}}}

\newcommand{\catt}{\textnormal{\textsf{catt}}\xspace}
\newcommand{\coh}[4][]{\ensuremath{\operatorname{coh}^{#1}(#2,#3,#4)}}
\DeclareMathOperator{\coherence}{coh}
\DeclareMathOperator{\var}{var}

\newcommand{\compcell}[3]{\ensuremath{#2\ast_{#1}#3}}
\DeclareMathOperator{\unbcomp}{comp}
\DeclareMathOperator{\eh}{eh}

\title{Hom \(\omega\)\-categories of a computad are free}

\author{ Thibaut
  Benjamin\footnote{University of Cambridge, \texttt{tjb201@cam.ac.uk}} \and
  Ioannis Markakis\footnote{University of Cambridge,
    \texttt{ioannis.markakis@cl.cam.ac.uk}}}

\begin{document}

\maketitle

\begin{abstract}
  We provide a new description of the hom functor on weak $\omega$\-categories,
  and we show that it admits a left adjoint that we call the suspension functor.
  We then show that the hom functor preserves the property of being free on a
  computad, in contrast to the hom functor for strict \(\omega\)\-categories.
  Using the same technique, we define the opposite of an
  $\omega$\-category with respect to a set of dimensions, and we show that this
  construction also preserves the property of being free on a computad.
  Finally, we show that the constructions of opposites and homs commute.
\end{abstract}

\input{1-intro}
\input{2-pds}
\input{3-categories}
\input{5-homs}
\input{4-opposites}
\input{6-applications}

\section*{Acknowledgements}

We would like to thank Prof. Jamie Vicary for his support during this project.
We would also like to thank Prof. John Bourke and Dimitri Ara for their helpful
remarks on earlier versions of this work. Furthermore, the second author would
like to acknowledge funding from the Onassis foundation - Scholarship ID: F ZQ
039-1/2020-2021.

\bibliographystyle{plainurl}
\bibliography{bibliography}

\end{document}

%% file: 1-intro.tex

\section{Introduction}

In recent years, higher category theory has found a wide range of applications
in various fields of mathematics and computer science. Globular higher
categories and their computads have been used significantly in rewriting
theory~\cite{mimram_3dimensional_2014}, homology
theory~\cite{lafont_polygraphic_2009}, homotopy
theory~\cite{garner_homomorphisms_2010} and topological quantum field
theory~\cite{bartlett_modular_2015}. Furthermore, globular higher groupoids
describe identity types in homotopy type theory~\cite{lumsdaine_weak_2009,
  altenkirch_syntactical_2012, vandenberg_types_2011} and they are conjectured
to model spaces~\cite{grothendieck_pursuing_1983,henry_homotopy_2023}. Computads
play an important role in the homotopy theory of higher categories. They are
combinatorial structures that freely generate higher categories, and
characterise cofibrant objects in the folk model structure of strict
\(\omega\)\-categories~\cite{metayer_cofibrant_2008,lafont_folk_2010}.

Strict \(\omega\)\-categories are one of the most studied kind of globular
higher category. They are usually defined via enrichment, so their composition
operations are associative and unital and satisfy the interchange laws. An
immediate consequence of their definition is the existence of opposite and hom
\(\omega\)\-categories. Opposites are obtained by exchanging the source and
target of cells and reversing the order of composition, whereas enrichment gives
an \(\omega\)\-category structure to the collection of cells between two
objects. While opposites preserve freeness on a computad, the construction of
hom \(\omega\)\-categories does not. This is due to the Eckmann-Hilton argument,
which imposes a commutativity relation on the arrows of certain hom
\(\omega\)\-categories. In this paper, we focus on weak \(\omega\)\-categories
and show that the construction of homs does preserve freeness on a computad in
this case.

\paragraph{Contributions.}
We give a novel description of the hom functor for weak \(\omega\)\-categories,
using the computads-first approach of Dean et al.~\cite{dean_computads_2024},
and we show that it admits a left adjoint, the suspension functor. Moreover,
we use the same approach to construct the opposites of \(\omega\)\-categories,
and we show that the two constructions commute.

The opposite and the suspension functors preserve freeness on a computads by
construction. Our main result is that the same is true for the hom functor. We
first identify the generators of the hom \(\omega\)\-category of a computad as
its indecomposable cells. Every cell of the hom \(\omega\)\-category can be
decomposed into indecomposable cells. Moreover, the generators satisfy no
relations, so the hom \(\omega\)\-category of a computad is freely generated by
by its indecomposable cells.

We further provide an implementation of opposites as a meta\-operation in the
proof assistant \catt\footnote{Available at
\url{https://github.com/thibautbenjamin/catt}} for \(\omega\)\-categories,
based on the type theory of Finster and
Mimram~\cite{finster_typetheoretical_2017}. This implementation relies on the
equivalence of contexts in this type theory with finite
computads~\cite{benjamin_catt_2024}. In conjunction with the automatic
computation of suspensions~\cite{benjamin_type_2020} and of
inverses~\cite{benjamin_invertible_2024}, this meta\-operation can
be used, for example, to reduce the size of the code that produces cells
witnessing the Eckmann\-Hilton argument, its inverse and their cancellation
witnesses by 95\%.

\paragraph{Related work.}
Both our constructions proceed similarly, taking advantage of the inductive
description of computads of Dean et al.~\cite{dean_computads_2024} to lift
adjoints first from globular sets to computads, and then from computads to
\(\omega\)\-categories. Those constructions provide further evidence that this
computads-first approach is well\-suited for developing the theory of
\(\omega\)\-categories.

The construction of the suspension functor is inspired by the suspension
meta\-operation on the type theory \catt, introduced in the first author's
thesis~\cite{benjamin_type_2020}. The former extends the latter via the the
identification of \catt contexts with finite
computads~\cite{benjamin_catt_2024}.

Hom \(\omega\)\-categories have been previously constructed by Cottrell and
Fujii~\cite{cottrell_hom_2022}, using Leinster's approach to
\(\omega\)\-categories. We give an alternative explicit construction of the hom
\(\omega\)\-category functor, that allows us to prove that it preserves the
property of being free on a computad. We note that even though both
constructions provide an \(\omega\)\-category structure on the same globular set
of morphisms, they need not coincide. We conjecture that instead they are merely
weakly equivalent.

\paragraph{Overview.}
We start in Section~\ref{sec:glob-past-diag} by recalling the suspension and
hom adjunction
\[ \susp : \glob \rightleftarrows \bipointed : \deloop \]
between globular sets and bipoitned globular sets. We then use this adjunction
to describe globular pasting diagrams, the arities of the opertions of
\(\omega\)\-categories. In Section~\ref{sec:comp}, we introduce computads and
weak \(\omega\)\-categories following Dean et al.~\cite{dean_computads_2024}.
First, the category of computads and morphisms of free \(\omega\)\-categories is
defined by an inductive construction, together with an adjunction with globular
sets
\[\free : \glob \rightleftarrows \comp : \cell.\]
Then, \(\omega\)\-categories are defined as algebras for the monad \(T\) induced
by this adjunction.

Section~\ref{sec:hom} is dedicated to the construction of the suspension and hom
adjunction. We start by extending the suspension of globular sets to computads
in a way that is compatible with the adjunction \(\free\dashv\cell\).
Compatibility with the adjunction gives rise to a natural transformation
\[\susp^T : \susp T \Rightarrow \Tbip \susp \]
where \(\Tbip\) is the free bipointed \(\omega\)\-category monad. The mate of
this natural transformation is a morphism of monads. Hence, it induces an
extension
\[
  \deloop : \catbip \to \omega\Cat
\]
of the hom functor to \(\omega\)\-categories by the formal theory of
monads~\cite{street_formal_1972}. Moreover, by the adjoint lifting theorem, we
obtain a left adjoint to this functor.

We then prove our main result that the hom functor preserves freeness on a
computad. From an \(\omega\)\-category \(\wcat X\), we construct a computad
\(C_{\wcat X}\) whose generators are the indecomposable cells of \(\wcat X\)
together with a moprhism of \(\omega\)\-categories
\(\sigma_{\wcat X} : C_{\wcat X} \to \wcat X\) from the free
\(\omega\)\-category on \(C_{\wcat X}\) to \(\wcat X\). We then show that
\(\sigma_{\wcat X}\) is an isomorphism in the case that \(\wcat X\) is the hom
\(\omega\)\-category of a computad. It follows that the hom \(\omega\)\-category
functor restricts to a functor
\[\deloop : \compbip \to \comp  \]
on the level of computads.

In Section~\ref{sec:opp}, we use the same inductive technique to construct the
opposites of a weak \(\omega\)\-category. In general, $n$\-categories admit
$2^n-1$ opposites, obtained by reversing the direction of cells of certain
dimensions. Similarly, $\omega$\-categories have uncountably many opposites,
indexed by the group \(G\) of subsets of the positive natural numbers. We
construct opposite \(\omega\)\-categories by first defining the opposites of
globular sets and computads
\begin{align*}
  \op : \glob &\to \glob &
  \op : \comp &\to \comp
\end{align*}
in a way that is compatible with the adjunction $\free\dashv\cell$. This gives
rise to an endomorphism of the free \(\omega\)\-category monad, hence an
endofunctor
\[\op : \omega\Cat\to \omega\Cat\]
of the category of $\omega$\-categories. We then show that the formation of
opposites gives rise to an action \(G\) on the category of
\(\omega\)\-categories, so in particular the functors \(\op\) are involutive.
Finally, we show that the formation of opposites and homs commute by proving a
commutativity result for the associated morphisms of monads.

Throughout this paper we illustrate with simple examples how the constructions
we introduce interract with the operations of \(\omega\)\-categories. In
Section~\ref{sec:eh}, we conclude with a discussion of a more involved example,
that of the Eckmann-Hilton cells. We explain how the suspension can be used to
raise the dimension of Eckmann-Hilton cells. We remark that certain opposites of
Eckmann-Hilton cells are also their inverses, which is proven in our subsequent
article~\cite{benjamin_invertible_2024}.


%% file: 2-pds.tex

\section{Globular pasting diagrams}
\label{sec:glob-past-diag}

In this section, we briefly recall the notion of globular pasting diagrams,
since they are a basic ingedient for any definition of weak
$\omega$\-categories. Those are a family of globular sets such that diagrams
indexed by them in a strict $\omega$\-category can be composed in a unique way.
Pasting diagrams are parametrised by rooted, planar
trees~\cite{batanin_monoidal_1998}, an inductive description of whose as
iterated lists was recently given by Dean et
al~\cite[Section~2]{dean_computads_2024}. Our presentation in this section
follows ibid. and Leinster~\cite[Appendix~F.2]{leinster_higher_2004}, noting
that pasting diagrams are bipointed globular sets generated by the suspension
and the wedge sum operations.

To set the notation, we recall that \emph{globular sets} are presheaves on the
category $\G$ of globes with objects the natural numbers $\N$ and morphisms
freely generated by the source and target inclusions
\[s_n,t_n\colon[n]\to[n+1]\]
under the globularity relations:
\begin{align*}
  s_{n+1}\circ s_n &= t_{n+1} \circ s_n & s_{n+1}\circ t_n &= t_{n+1} \circ t_n.
\end{align*}
In other words, a globular set $X$ consists of a set $X_n$ for every natural
number $n\in \N$ together with \emph{source} and \emph{target} functions
\[\src,\tgt \colon X_{n+1}\to X_n\]
satisfying the duals relations:
\begin{align*}
  \src\circ\src &= \src\circ \tgt & \tgt\circ\src &=\tgt\circ\tgt.
\end{align*}
We will call elements of $X_n$ the \emph{$n$\-cells of $X$}. The
\emph{$k$\-source} and \emph{$k$\-target} of an $n$\-cell $x$ for $k < n$ are
the $k$\-cells defined by
\begin{align*}
    \src_k x &= \src(\cdots(\src x)) & \tgt_k x &= \tgt(\cdots(\tgt x))
\end{align*}
We will denote by $\disk n$ the representable globular set associated to a
natural number $n\in \N$, and call it the \emph{$n$\-disk}.

Bipointed globular sets are triples $(X,x_-,x_+)$ consisting of a
globular set and two distinguished $0$\-cells $x_-,x_+$ of it. They form a
category $\bipointed$ toether with morphisms of globular sets that preserve the
distinguished $0$\-cells. By the Yoneda lemma, this is the coslice category
$\coslice \glob {\disk 0 + \disk 0}$ or the category of cospans of globular sets
from $\disk 0$ to itself.

The category of bipointed globular sets is locally finitely presentable as a
coslice of a presheaf topos \cite[Proposition~1.57]{adamek_locally_1994}, so in
particular it is complete and cocomplete. Limits and connected colimits in
$\bipointed$ are computed as in $\glob$, i.e. they are created by the
functor $\bipointed \to \glob$ that forgets the basepoints. The
coproduct of a family of bipointed sets is computed as the wide pushout of the
corresponding maps out of $\disk 0 + \disk 0$.

Being a category of cospans, the category $\bipointed$ is
monoidal with respect to the composition of cospans, which we will call the
\emph{wedge sum}. More explicetely, the wedge sum $\vee$ of a pair of bipointed
globular sets $(X,x_-,x_+)$ and $(Y, y_-,y_+)$ is obtained by the following
pushout square in $\glob$,
\[
  \begin{tikzcd}[column sep = small, row sep = small]
    & & X\vee Y & & \\
    & X
    \ar["\inc_1", ru, dashed]
    & &
    Y\ar["\inc_2"', lu, dashed] & \\
    \disk{0}
    \ar[ru,"x_-"] & &
    \disk{0}
    \ar[ru,"y_-"]
    \ar[lu,"x_+"']
    \ar[uu,phantom,"\llcorner"{rotate = 45, very near end}] & &
    \disk{0}
    \ar[lu,"y_+"']
  \end{tikzcd}
\]
with basepoints the image of $x_-$ and the image of $y_+$ in the pushout. The
unit of the wedge sum is given by the $0$\-disk $\disk 0$ with both
basepoints being its unique $0$\-cell. More generally, we will denote by
$\bigvee_{i=1}^n X_i$ the iterated monoidal product of a finite family of
bipointed globular sets $X_1, \ldots, X_n$, and we will denote the inclusion of
the $j$\-th component for $1\le j\le n$ by
\[\inc_j \colon X_j \to \bigvee_{i=1}^n X_j. \]
This is a morphism of globular sets, that is only a morphism of bipointed
globular sets for $n = 1$.

The \emph{suspension} of a globular set $X$ is the bipointed globular set
$\susp X$ with two $0$\-cells $v_-$ and $v_+$, and with positive dimensional
cells given by
\[ (\susp X)_{n+1} = X_n \]
for every $n\in \N$. The source and target maps of an $n$\-cell of $\susp X$
for $n > 2$ is given by its source and target in $X$, while the source and
target of $1$\-cells are given by $v_-$ and $v_+$ respectively. The basepoints
of the suspensions are $v_-$ and $v_+$. Suspension is left adjoint to the
\emph{path space} functor $\deloop \colon \bipointed \to \glob$ sending a bipointed
globular set $(X,x_-,x_+)$ to the globular set given by
\[
  \deloop (X,x_-,x_+)_n =
    \set{ x\in X_{n+1} \such \src_0(x) = x_- \text{ and } \tgt_0(x) = x_+ }
\]
The unit of the adjunction is the identity of the functor
\[\deloop \susp = \id,\]
while the counit $\kappa : \susp\deloop \Rightarrow \id$ is the natural transformation
with components the bipointed morphisms
\[\kappa : \susp\deloop (X,x_-,x_+) \to (X, x^-, x^+)\]
given by the subset inclusions $\deloop (X,x_-,x_+)_n \subseteq X_{n+1}$.

Finally, we have introduced all the ingredients to define globular pasting
diagrams and the family parametrising them. We will call elements of that
family Batanin trees following Dean et al~\cite{dean_computads_2024}.

\begin{definition}
    A \emph{Batanin tree} is a list
  \(\bt{B_{1},\ldots,B_{n}}\), where the \(B_{i}\) are Batanin trees.
\end{definition}

In other words, the set $\bat$ of Batanin trees is the carrier of the initial
algebra of the list endofunctor $\List \colon \Set \to \Set$ given by
\[
  \List X = \coprod_{n\in \N} X^n
\]
with the obvious action on morphisms. In particular, there exists a tree $\bt{}$
corresponding to the empty list, and using this tree, we can define more
complicated trees, such as the tree \[B=\bt{\bt{\bt{},\bt{}},\bt{}}.\]
It is convenient to visualise Batanin trees as planar trees by representing
$\bt{}$ as a tree with one root and no branches, and $\bt{B_1, \ldots, B_n}$ as
a tree with a new root and $n$ branches, each of which is connected to the
root of the tree corresponding to $B_i$. For example, the tree $B$ above can be
visualised as
\[
  \begin{tikzcd}[column sep = small]
    \bullet & & \bullet\\
    & \bullet \ar[ul,-]\ar[ur,-] & & \bullet \\
    & & \bullet \ar[ul,-] \ar[ur,-]
  \end{tikzcd}
\]

The \emph{dimension} of a Batanin tree is the height of the corresponding planar
tree, or equivalently the maximum of the dimension of the cells in the
corresponding globular pasting diagram, defined below. It can be computed
recursively by
\[
  \dim (\bt{B_1, \ldots, B_n}) = \max(\dim B_1 +1, \ldots, \dim B_n +1).
\]
In particular, it follows that $\bt{}$ is the unique tree of dimension $0$.

\begin{defi}\label{def:pos}
  The \emph{bipointed globular set of positions} of a Batanin tree $B$ is the
  bipointed globular set $\Pos B$ defined recursively by
  \[
    \Pos{\bt{B_{1},\ldots,B_{n}}} =
    \bigvee_{i=1}^{n} \susp\Pos {B_{i}}.
  \]
  The \emph{globular pasting diagram} of a Batanin tree $B$ is the
  underlying globular set of $\Pos B$, according to the fomulae
  in~\cite[Appendix~F.2]{leinster_higher_2004}.
\end{defi}

A way to calculate the globular set of positions of a tree is described
in~\cite{berger_cellular_2002}, where positions correspond to \emph{sectors} of
the tree, i.e. the spaces between two consecutive branches at each node, as well
as the space before the first branch and the one after the last one. Under this
description, the basepoint are given by the left-most and right-most sector at
the root. For the tree $B$ above, we can label the position as follows.
\[
  \begin{tikzcd}[column sep = small]
    \phantom{\bullet} & \overset{a}{\bullet}
    & \phantom{\bullet} & \overset{b}{\bullet} & \phantom{\bullet}\\
    & \phantom{\bullet} & \bullet \ar[ul,-]\ar[ur,-]
    \ar[ull,phantom,"{\scriptstyle f}"very near start]
    \ar[u,phantom,"{\scriptstyle g}"very near start]
    \ar[urr,phantom,"{\scriptstyle h}"very near start] &
    \phantom{\bullet}
    & \bullet \ar[u,phantom,"{\scriptstyle k}"very near start]
    & \phantom{\bullet} \\
    & & & \bullet \ar[ul,-] \ar[ur,-] \ar[ull,phantom,"{\scriptstyle
      x}"very near start] \ar[u,phantom,"{\scriptstyle y}"very near
    start] \ar[urr,phantom,"{\scriptstyle z}"very near start]
  \end{tikzcd}
\]
The dimension of a position is given by the distance of the node it is attached
in from the root, while its source and target are given by the positions right
below it. Therefore, the globular set of positions of $B$ is the following
globular set
\[
  \begin{tikzcd}
    x\ar[r,"g"{description}]\ar[r,bend left = 50,"f"]\ar[r,bend right
    = 50,"h"below] \ar[r,bend left = 25, phantom, "\Downarrow{\scriptstyle
    a}"] \ar[r,bend right = 25, phantom, "\Downarrow{\scriptstyle
    b}"] & y\ar[r,"k"] & z
  \end{tikzcd}
\]
which is bipointed by the positions $x$ and $z$ respectively. Here, the
positions $f,g,h,a,b$ are the positions of the left branch of $B$, while $k$ is
the position of its right branch. The dimension of those positions has been
raised by the suspension operation. The $0$\-positions $x,y,z$ are the new cells
created by the suspension operation. The two basepoints of $\Pos{B}$ are
given by $x$ and $z$.

\begin{definition}
  The \emph{$k$\-boundary} of a Batanin tree $B$ is the tree $\partial_k B$
  defined recursively by
  \begin{align*}
    \partial_0 B
      &= \bt{} \\
    \partial_{k+1} \bt{B_1,\dots,B_n}
      &= \bt{\partial_k B_1,\dots,\partial_k B_n}
  \end{align*}
\end{definition}

\noindent
The \emph{$k$\-boundary} of a tree $B$ is the tree obrained by removing all
nodes of $B$ whose distance from the root is at least $k$. In terms of
pasting diagrams, this amounts to removing all cells of dimension more than $k$
and identifying all parallel $k$\-cells. For example, the $1$\-boundary of the
tree $B$ considered above is the following tree.
\[
  \bdry[1]{B} =
  \begin{tikzcd}[column sep = small]
    & \bullet & & \bullet \\
    & & \bullet \ar[ul,-] \ar[ur,-]
  \end{tikzcd}
  \qquad\qquad
  \Pos{\bdry[1]{B}} =
  \begin{tikzcd}
    \bullet \ar[r] & \bullet\ar[r] & \bullet
  \end{tikzcd}
\]

The positions of the boundary can be included back into the positions of the
original tree in two ways, the \emph{source} and \emph{target inclusions}
\[
  \srcps[k]{B}, \tgtps[k]{B} \colon \Pos{\bdry[k]{B}}\to \Pos{B}
\]
defined recursively as follows: the morphisms $\srcps[0]{B}$ and $\tgtps[0]{B}$
out of $\disk 0 \cong \Pos{\bt{}}$ select the first and second basepoint
respectively, while for $B = \bt{B_1,\dots,B_n}$ the morphisms
$\srcps[k+1]{B}$ and $\tgtps[k+1]{B}$ are given by
\begin{align*}
  \srcps[k+1]{B} &= \bigvee_{i=1}^n \susp \srcps[k]{B_i} &
  \tgtps[k+1]{B} &= \bigvee_{i=1}^n \susp \tgtps[k]{B_i}.
\end{align*}
In particular, the source and target inclusions are morphisms of bipointed
globular sets when $k>0$.

As an illustration for Batanin tree and pasting schemes, we consider a few
important families of Batanin trees, for which we also illustrate a few of the
corresponding globular pasting diagrams for visualisation:
\begin{example}\label{ex:ps-disk}
  We define the \emph{\(k\)\-dimensional disk} tree \(D_{k}\) recursively on
  \(k\) by
    \begin{align*}
      D_{0} &= \bt{} & D_{k+1} = \bt{D_{k}}.
    \end{align*}
  A simple inductive argument shows that \(\Pos{D_{k}} \cong \disk k\).
\end{example}

\begin{example}\label{ex:ps-comp}
  We now define a family of Batanin trees that are the arities of the binary
  compositions in \(\omega\)\-categories.
  We define first for \(n,m > 0\),
  \begin{align*}
    B_{n,0,m} &= \bt{D_{n},D_{m}} & \Pos{B_{n,0,m}} &\cong \disk n \vee \disk m.
  \end{align*}
  A few low\-dimensional globular pasting diagrams in this family are
  illustrated below
    \[
      \begin{array}{wc{.3\textwidth}wc{.3\textwidth}wc{.3\textwidth}}
        \Pos{B_{1,0,1}}
        & \Pos{B_{2,0,1}}
        & \Pos{B_{3,0,3}} \\
        \begin{tikzcd}[ampersand replacement=\&]
          \bullet\ar[r] \& \bullet\ar[r] \& \bullet
        \end{tikzcd}
        &
          \begin{tikzcd}[ampersand replacement=\&]
            \bullet
            \ar[r, bend left]
            \ar[r, bend right]
            \ar[r,phantom,"\Downarrow"]
            \& \bullet
            \ar[r]
            \& \bullet
          \end{tikzcd}
        &
          \begin{tikzcd}[ampersand replacement=\&]
            \bullet
            \ar[r, bend left]
            \ar[r, bend right]
            \ar[r,phantom,"\Downarrow\Rrightarrow\Downarrow"]
            \& \bullet
            \ar[r, bend left]
            \ar[r, bend right]
            \ar[r,phantom,"\Downarrow\Rrightarrow\Downarrow"]
            \& \bullet
          \end{tikzcd}
      \end{array}
    \]
    We then define recursively on \(k\) the family
    \[B_{n,k,m} = \bt{B_{n-1,k-1,m-1}}\]
    for every \(n,m > 0\) and \(0\leq k < \min(n,m)\). Using that the suspension
    preserves representables and connected colimits, we conclude that
    \[
      \Pos {B_{n,k,m}} \cong \disk n \amalg_{\disk k} \disk m
    \]
    where the pushout is over the target and source inclusions respectively. We
    illustrate some of the globular pasting diagrams in this family below.
        \[
      \begin{array}{wc{.3\textwidth}wc{.3\textwidth}wc{.3\textwidth}}
        \Pos{B_{2,1,2}}
        & \Pos{B_{3,1,2}}
        & \Pos{B_{3,2,3}} \\
        \begin{tikzcd}[ampersand replacement=\&]
          \bullet
          \ar[rr, bend left=50]
          \ar[rr,phantom,"\Downarrow", bend left=25]
          \ar[rr]
          \ar[rr,phantom,"\Downarrow", bend right=25]
          \ar[rr, bend right=50]
          \&\& \bullet
        \end{tikzcd}
        &
          \begin{tikzcd}[ampersand replacement=\&]
            \bullet
            \ar[rr, bend left=50]
            \ar[rr,phantom,"\Downarrow\Rrightarrow\Downarrow", bend left=25]
            \ar[rr]
            \ar[rr,phantom,"\Downarrow", bend right=25]
            \ar[rr, bend right=50]
            \&\& \bullet
          \end{tikzcd}
        &
          \begin{tikzcd}[ampersand replacement=\&]
            \bullet
            \ar[rr, bend left]
            \ar[rr, phantom,"\Downarrow\Rrightarrow\Downarrow\Rrightarrow\Downarrow"]
            \ar[rr, bend right]
            \&\& \bullet
          \end{tikzcd}
      \end{array}
    \]
    A simple induction shows that the dimension of the Batanin tree
    \(B_{n,k,m}\) is given by
    \[
      \dim(B_{n,k,m}) = \max(n,m)
    \]
\end{example}


%% file: 3-categories.tex

\section{Computads and \texorpdfstring{$\omega$}{ω}-categories}
\label{sec:comp}

Dean et al.~\cite{dean_computads_2024} recently presented a new definition of
$\omega$\-categories and their computads, inspired by the type\-theoretic
definition of Finster and Mimram~\cite{finster_typetheoretical_2017}, and they
showed that their notion of $\omega$\-category coincides with the operadic
definition of Leinster~\cite{leinster_higher_2004}. In this approach, first
a category of computads $\comp$ is defined together with an adjunction
\[
  \free \colon \glob \rightleftarrows \comp \colon \cell
\]
and then $\omega$\-categories are defined as algebras for the monad
\[
  T \colon \glob \to \glob
\]
induced by the adjunction. We recall that morphisms of computads here
are strict $\omega$\-functors, and not Batanin's morphisms of
computads~\cite{batanin_computads_1998}. In other words, the
comparison functor
\[
  K^T \colon \comp \to \omega\Cat
\]
is fully faithful and injective on objects.

We will briefly recall the definition of computads and the $\free\dashv\cell$
adjunction. First, categories $\comp_n$ of $n$\-computads are defined
recursively for every natural number $n\in \N$, together with forgetful functor
\[ u_n \colon \comp_n \to \comp_{n-1}\]
for $n>0$. In the same mutual recursion, functors
\begin{align*}
  \free_n &\colon \glob \to \comp_n \\
  \cell_n &\colon \comp_n \to \Set \\
  \type_n &\colon \comp_n \to \Set
\end{align*}
are defined and natural transformations
\begin{align*}
  \ty_n &\colon \cell_n \Rightarrow \type_nu_n \\
  \pr_i &\colon \type_n \Rightarrow \cell_n
\end{align*}
for $i = 1, 2$. Here the functors $\cell_n$ and $\type_n$ return the set of
$n$\-cells, and the set of pairs of parallel $n$\-cells of the
$\omega$\-category generated by a computad $C$, while $\ty_n$ returns the source
and the target of an $n$\-cell.

An \emph{$n$\-computad} is a triple $C$ consisting of an $(n-1)$\-computad
$C_{n-1}$, a set of $n$\-dimensional \emph{generators} $V_n^C$ and an attaching
function
\[ \phi_n^C \colon V_n^C \to \type_{n-1}(C_{n-1}) \] assigning to each generator
a source and target. A morphism $\sigma \colon C\to D$ consists of a morphism
$\sigma_{n-1} \colon C_{n-1} \to D_{n-1}$ and a function
$\sigma_{n,V} \colon V_n^C\to \cell_n D$ compatible with the source and target
functions in the sense defined in~\cite[Section~3.1]{dean_computads_2024}. The
forgetful functors $u_n$ are the obvious projections. As a base case for this
definition, here we let $\comp_{-1}$ be the terminal category and $\type_{-1}$
the functor choosing some terminal set.

The set $\cell_n C$ of \emph{$n$\-cells of a computad} $C$ is inductively
defined together with the set of morphisms with target $C$ and the function
$\ty_{n,C}$. Cells of $C$ are either of the form $\var v$ for a generator
$v\in V_n^C$, or when $n>0$, they are \emph{coherence cells} $\coh{B}{A}{\tau}$,
where $B$ is a tree of dimension at most $n$, $A$ is an $(n-1)$\-sphere of
$\free_{n-1} \Pos B$, satisfying a \emph{fullness} condition that will be
explained below, and $\tau \colon \free_n\Pos{B} \to C$ is a morphism. The
boundary of a cell is given recursively by the formula
\begin{align*}
  \ty_{n,C}(\var v) &= \phi_n^C(v) \\
  \ty_{n,C}(\coh{B}{A}{\tau}) &= \type_{n-1}(\tau_{n-1})(A)
\end{align*}
The functor $\free_n$ sends a globular set $X$ to the computad
\begin{align*}
  \free_{n}X&=(\free_{n-1}X, X_{n}, \phi_{n}^{X}) &
  \phi_n^X(x) &= (\var (\src x), \var (\tgt x)),
\end{align*}
and a morphism $f \colon X\to Y$ to the morphism consisting of $\free_{n-1}f$
and $\var\circ f_n$.

The functor $\type_n$ sends an $n$\-computad $C$ to the set
\[
  \type_n C =
  \set{ (a,b) \in \cell_n C\times \cell_n C \such \ty_na = \ty_n b}
\]
and acts on morphisms in the obvious way. The projection natural transformations
are the obvious ones. We will denote by
\[\src, \tgt \colon \cell_n \Rightarrow \cell_{n-1}u_n\]
the composite of $\ty_n$ with the projections.

The fullness condition mentioned above for
$A = (a,b) \in \type_n\free_n\Pos{B}$ is a condition on the generators used to
define $a$ and $b$. It is equivalent to the statement that
\begin{align*}
  a &= \cell_n\free_n(s_n^B)(a') &
  b &= \cell_n\free_n(t_n^B)(b')
\end{align*}
for cells $a',b'$ of $\free_{n} \Pos{\partial_n B}$ using all generators of
$\partial_n B$. That means that the \emph{support} of $a',b'$
contains all positions of $\partial_n B$, where the support of an $n$\-cell $c$
over a computad $C$ is the set of generators defined by
\begin{align*}
  \supp(\var v) &= \begin{cases}
    \{v\},
      &\text{when }n=0 \\
    \{v\}\cup\supp(\pr_1\phi_n^C(v))\cup\supp(\pr_2\phi_n^C(v)),
      &\text{when }n>0
  \end{cases} \\
  \supp(\coh{B}{A}{\tau})
    &= \bigcup_{k\le n}\bigcup_{v\in \Pos[k]{B}} \supp(\tau_{k,V}(v))
\end{align*}

This completes the inductive definition. The category $\comp$ of computads is
the limit of the categories $\comp_n$ for all $n\in\N$, i.e computads
$C = (C_n)_{n\in \N}$ are sequences of $n$\-computads $C_n$ such that
$u_{n+1}C_{n+1} = C_n$, and moprhisms of such are sequences of morphisms. The
free functor
\[
  \free \colon \glob \to \comp
\]
is the functor with components $\free_n$ for all $n\in \N$, while the cell
functor
\[
  \cell \colon \comp\to \glob
\]
sends a computad $C$ to the globular set consisting of $\cell_n C_n$ for all
$n\in \N$, and the source and target functions defined above. The unit and the
counit of the adjunction
\begin{align*}
  \eta &\colon \id \Rightarrow \cell\free \\
  \varepsilon &\colon \free\cell \Rightarrow \id
\end{align*}
are described as follows: The unit \(\eta\) sends a cell $x$ of a globular set
$X$ to the generator $\var x$, the counit \(\varepsilon\) consists of the
morphisms $\varepsilon_C \colon \free\cell C\to C$ given for all \(n\in\N\) by
the identities of the set
\[V_n^{\free\cell C} = \cell_n C. \]

By definition, $\omega$\-categories are algebras for the monad $(T,\mu,\eta)$ on
$\glob$ induced by the adjunction $\free\dashv \cell$. In particular, there
exists a free/underlying adjunction
\[
  F^T : \glob \rightleftarrows \omega\Cat : U^T
\]
between globular sets and $\omega$\-categories, and there exists a comparison
functor
\[  K^T : \comp \to \omega\Cat \]
sending a computad $C$ to the $\omega$\-category
$(\cell C, \cell \varepsilon_C)$. Moreover, $K^T$ is a morphism of adjunctions
meaning that
\begin{align*}
  F^T &= K^T \free & \cell &= U^T K^T,
\end{align*}
and it is fully faithful~\cite[Proposition~4.6]{dean_computads_2024}.

\subsection{Binary compositions in \texorpdfstring{\(\omega\)}{ω}-categories}
\label{sec:cell-compositions}
The aim of this section is to illustrate this definition by describing the
binary compositions of \(\omega\)\-categories. To achieve this, we construct a
cell \(\unbcomp_{n,k,m}\) of dimension \(\max(n,m)\) in the free
\(\omega\)\-category on the globular pasting diagram \(\Pos{B_{n,k,m}}\),
introduced in Example~\ref{ex:ps-comp}. By the pushout formulae in
Example~\ref{ex:ps-comp}, morphisms of \(\omega\)\-categories
\(f : T\Pos{B_{n,k,m}} \to \wcat{X}\) correspond precisely to pairs of
an \(n\)\-cell \(c\) and an \(m\)\-cell \(c'\) in \(\wcat X\) such that
\(\tgt_{k}(c) = \src_{k}(c')\). The \(k\)\-composite of \(c\) and \(c'\) is then
defined to be the cell
\[
  \compcell{k}{c}{c'} = f(\unbcomp_{n,k,m}).
\]
The rest of the operations of \(\omega\)\-categories can be obtained in a
similar way, by defining specific cells in well-chosen computads. This is in
particular the case for associators, unitors and interchangers. We now define
the cells \(\unbcomp_{n,k,m}\) starting from simpler cases and building up
towards the general case.

We first define the cell \(\unbcomp_{1,0,1}\). By the pushout formula, we
visualise the globular pasting diagram \(\Pos{B_{1,0,1}}\) as follows:
\[
  \begin{tikzcd}
    x\ar[r,"f"] & y\ar[r,"g"] & z
  \end{tikzcd},
\]
we then define the \(1\)\-cell of \(\free\Pos{B_{1,0,1}}\)
\[
  \unbcomp_{1,0,1} :=\coh{B_{1,0,1}}{(\var x,\var z)}{\id},
\]
with source \(\var x\) and target \(\var z\). Similarly, we use that
\(\bdry[k] B_{k+1,k,k+1} = D_k\) to define the cell
\[
  \unbcomp_{k+1,k,k+1} :=\coh{B_{k+1,k,k+1}}{(\var x_{k+1},\var z_{k+1})}{\id},
\]
where \(x_{k+1}\) and \(z_{k+1}\) are respectively the images of the
top\-dimensional cell of the disk \(\Pos{D_{k}}\) under the source and target
inclusions \(\Pos{D_{k}} \to \Pos{B_{k+1,k,k+1}}\) defined in
Section~\ref{sec:glob-past-diag}. The source and target of
\(\unbcomp_{k+1,k,k+1}\) imply that the source of \(\compcell{k}{c}{c'}\) is
\(\src(c)\) and its target is \(\tgt(c')\), as expected of composition
operations.

We then proceed to define the \(n\)\-cell \(\unbcomp_{n,k,n}\) recursively on
\(n-k\). The base case \(n = k+1\) has been defined above. For \(n>k+1\), we
have that \(\bdry[n-1] B_{n,k,n} = B_{n-1,k,n-1}\). We then define the cell
\(\unbcomp_{n,k,n}\) to be
\[
  \coh{B_{n,k,n}}{(T(s_{n-1}^{B_{n,k,n}})(\unbcomp_{n-1,k,n-1}),
    T(t_{n-1}^{B_{n,k,n}})(\unbcomp_{n-1,k,n-1}))}{\id}.
\]
It follows from this definition that the source and target of the
\(k\)\-composite \(\compcell{k}{c}{c'}\) are respectively the \(k\)\-composite
of the sources and the \(k\)\-composite of the targets.

Finally, we define the cell \(\unbcomp_{n,k,m}\) for arbitrary \(n\neq m\) and
\(k < \min(n,m)\). When \(n > m\), we have that
\(\bdry[n-1] B_{n,k,m} = B_{n-1,k,m}\), and we define recursively the
\(n\)\-cell \(\unbcomp_{n,k,m}\) to be
\[
  \coh{B_{n,k,m}}{(T(s_{n-1}^{B_{n,k,m}})(\unbcomp_{n-1,k,m}),
    T(t_{n-1}^{B_{n,k,m}})(\unbcomp_{n-1,k,m}))}{\id}.
\]
When \(n < m\), we construct the \(m\)\-cell \(\unbcomp_{n,k,m}\) in a similar
way. The binary \(k\)\-composite \(\compcell{k}{c}{c'}\) is the whiskering of an
\(n\)\-cell with an \(m\)\-cell. Its source and target are the \(k\)\-composites
of the source and target of the higher dimensional cell with the lower
dimensional one. In the following sections, we will see that leveraging the
suspension and opposite operations, the construction of the cells
\(\unbcomp_{n,k,m}\) can be significantly simplified.


%% file: 5-homs.tex

\section{The suspension and hom functors}
\label{sec:hom}
Strict $\omega$\-categories are precisely categories enriched over strict
$\omega$\-categories, so for every strict $\omega$\-category $X$ and every pair
of $0$\-cells $x_-,x_+\in X_0$, the globular set $\deloop (X,x_,x_+)$ of cells
from $x_-$ to $x_+$ admits an $\omega$\-category structure in a functorial
way~\cite{leinster_higher_2004}. The same result was recently proven for
arbitrary $\omega$\-categories~\cite{cottrell_hom_2022} using the operadic
definition of Leinster. We provide an elementary construction of hom
\(\omega\)\-categories, and show that it preserves the property of being free on
a computad. This diverges from strict \(\omega\)\-categories, where the
Eckmann-Hilton argument is an obstruction to freeness of hom
\(\omega\)\-categories.

Recall that the \emph{hom} functor $\deloop : \bipointed\to \glob$, taking a
bipointed globular set to the globular set of cells from the first base-point to
the second one admits a left adjoint, the \emph{suspension} functor,
$\susp : \glob \to \bipointed$ with unit the identity and counit
$\kappa : \susp\deloop \Rightarrow \id$ given by subset inclusions. Our goal in
this section will be to lift this adjunction to computads and to
$\omega$\-categories.

\subsection{The suspension of a computad}
Let $\compbip$ be the category of computads with two chosen $0$\-cells and
morphisms preserving those $0$\-cells. By the $\free\dashv \cell$ adjunction and
the Yoneda lemma, this is precisely the slice of $\comp$ under
$\free (\disk 0 + \disk 0)$. Moreover, the adjunction descends to the slices to
give an adjunction
\[
  \freebip : \bipointed \rightleftarrows \compbip : \cellbip
\]
where $\freebip (X, x_-, x_+)$ is the computad $\free X$ with the $0$\-cells
$\var x_-$ and $\var x_+$, and $\cellbip(C, c_-,c_+)$ is the globular set
$\cell C$ with the basepoints $c_-$ and $c_+$ respectively. We now define a
suspension functor
\[
  \susp : \comp \to \compbip
\]
which generalises the suspension of globular sets.

We define first the \emph{suspension} of a Batanin tree $B$ to be the Batanin
tree
\[ \susp B = \bt{B}\] since by definition
\[
  \Pos{\susp B} = \susp \Pos{B}.
\]
We will then proceed inductively on $n\in \N$ to define a functor and two
natural transformations
\begin{align*}
  \susp &: \comp_n \to \comp_{n+1} \\
  \susp^{\cell} &: \cell_n \Rightarrow \cell_{n+1}\susp \\
  \susp^{\type} &: \type_n \Rightarrow \type_{n+1}\susp
\end{align*}
satisfying the following properties:
\begin{enumerate}[label={\textsc{(s{\footnotesize\arabic*})}}]
\item\label{susp-1} the suspension commutes with the forgetful functors, and the
  inclusion of globular sets into computads:
  \[\begin{tikzcd}
    {\comp_{n+1}} & {\comp_{n+2}} \\
    {\comp_{n}} & {\comp_{n+1}}
                  \arrow["{u_{n+1}}", from=1-1, to=2-1]
                  \arrow["{u_{n+2}}", from=1-2, to=2-2]
                  \arrow["\susp", from=1-1, to=1-2]
                  \arrow["\susp"', from=2-1, to=2-2]
  \end{tikzcd}\qquad\qquad
  \begin{tikzcd}
    {\comp_n} & {\comp_{n+1}} \\
    \glob & \glob \arrow["{\free_n}", from=2-1, to=1-1] \arrow["{\free_{n+1}}"',
            from=2-2, to=1-2] \arrow["\susp"', from=2-1, to=2-2] \arrow["\susp",
            from=1-1, to=1-2]
  \end{tikzcd}\]
\item\label{susp-2} the natural transformations are compatible with the boundary
  natural transformations:
  \[\begin{tikzcd}[column sep = large]
    {\cell_{n+1}}
    && {\cell_{n+2}\susp} \\
    {\type_n u_{n+1}}
    & {\type_{n+1}\susp u_{n+1}}
    & {\type_{n+1} u_{n+1}\susp}
    \arrow["{\ty_{n+1}}"', Rightarrow, from=1-1, to=2-1]
    \arrow["{\susp^{\type}u_{n+1}}"', Rightarrow, from=2-1, to=2-2]
    \arrow["{\susp^{\cell}}", Rightarrow, from=1-1, to=1-3]
    \arrow["{\ty_{n+1}\susp}", Rightarrow, from=1-3, to=2-3]
    \arrow[equals, from=2-2, to=2-3]
  \end{tikzcd}\]
\item\label{susp-3} the natural transformations are compatible with the
  projection natural transformations for $i = 1, 2$:
  \[\begin{tikzcd}
    {\type_n} & {\type_{n+1}\Sigma} \\
    {\cell_n} & {\cell_{n+1}\Sigma}
                \arrow["{\pr_i\Sigma}", Rightarrow, from=1-2, to=2-2]
                \arrow["{\pr_i}", Rightarrow, from=1-1, to=2-1]
                \arrow["{\susp^{\cell}}"', Rightarrow, from=2-1, to=2-2]
                \arrow["{\susp^{\type}}", Rightarrow, from=1-1, to=1-2]
  \end{tikzcd}\]
\item\label{susp-4} the natural transformation $\susp^{\cell}$ preserves
  generators, in that for every globular set $X$ and $x\in X_n$, we have that
  \[
    \susp^{\cell}(\var x) = \var x
  \]
\item\label{susp-5} the natural transformation $\susp^{\type}$ preserves
  fullness, in that for every full $n$\-sphere $A$ of $\free\Pos{B}$, we have
  that $\susp^{\type}A$ is a full $(n+1)$\-sphere of $\free\Pos{\susp B}$.
\end{enumerate}

To start the induction, we recall that $\comp_{-1}$ is the terminal category and
that $\comp_0$ is the category $\Set$ of sets. The suspension functor is defined
as the functor picking the $2$\-element set $\set{v_-,v_+}$. We recall also that
the unique $(-1)$\-computad has a unique $(-1)$\-sphere. We define
$\susp^{\type}$ to be the natural transformation picking the $0$\-sphere
$(v_-,v_+)$. This concludes the base case. We will now assume that we have
defined the the data satisfying the properties we have cited, up to dimension
\(n-1\), for a fixed \(n\in\N\).

\paragraph{Computads.}
We will first define the functor \(\Sigma\) on all objects: Given an
$n$\-computad $C = (C_{n-1}, V_{n}^C, \phi_n^C)$, its suspension is the
$(n+1)$\-computad consisting of $\susp C_{n-1}$, the same set of generators, and
the attaching function $\phi_{n+1}^{\susp C}$ given by the composite
\[
  \phi_{n+1}^{\susp C} : V_n^C \xrightarrow{\phi_n^C} \type_{n-1}C_{n-1}
  \xrightarrow{\susp^{\type}} \type_n\susp C_{n-1}
\]
By definition, the suspension functor on objects commutes with the forgetful
functors. Using that the suspension functor on $(n-1)$\-computads commutes also
with the inclusion of globular sets, and that $\susp^{\cell}$ and hence
$\susp^{\type}$ preserve generators, we see that the suspension functor on
$n$\-computads also commutes with the inclusions.

\paragraph{Cells and morphisms.}
We then define the suspension of a morphism of \(n\)-computads together with the
natural transformation \(\susp^{\cell}\) mutually inductively, while showing
that Property~\ref{susp-2} holds. Given a computad \(C\), we define mutually
recursively \(\susp^{\cell}_{C}\) and the morphism
and \(\susp(\sigma)\) for every \(\sigma\) with target \(C\). For a generator
$v\in V_n^C$, we let
\[
  \susp^{\cell}_C(\var v) = \var v,
\]
and we compute that
\begin{align*}
  (\ty_{n+1,\susp C}\susp^{\cell}_C)(\var v)
  = \phi_{n+1}^{\susp C}(v)
  = (\susp^{\type}_{C_{n-1}}\ty_{n,C})(\var v).
\end{align*}
For a coherence $n$\-cell $c = \coh{B}{A}{\tau}$ of $C$, we let
\[
  \susp_C^{\cell}(\coh{B}{A}{\tau}) = \coh{\susp B}{\susp^{\type} A}{\susp \tau}
\]
using that the suspension commutes with $\free_n$ and $\Pos{-}$, and that it
preserves fullness. Then by naturality of $\susp^{\type}$, we compute that
\[
  (\ty_{n+1,\susp C}\susp^{\cell}_C)(\coh{B}{A}{\tau}) =
  (\susp^{\type}_{C_{n-1}}\ty_{n,C})(\coh{B}{A}{\tau}).
\]
For a morphism $\sigma : D\to C$, we let $\susp\sigma : \susp D\to \susp C$
consist of $\susp \sigma_{n-1}$ and the function
\[
  (\susp\sigma)_V : V_{n+1}^{\susp D} = V_n^D \xrightarrow{\sigma_V} \cell_n(C)
  \xrightarrow{\susp^{\cell}_C} \cell_{n+1}\susp C
\]
This is a well-defined morphism by the observation on the boundary of
$\susp^{\cell}$. Functoriality of the suspension, and naturality of
$\susp^{\cell}$ can be shown mutually inductively.

\paragraph{Spheres.}
Finally, the natural transformation $\susp^{\type}$ is defined for a computad
$C$ and an $n$\-sphere $ A = (a,b)$ of it again by
\[\susp^{\type}(a,b) = (\susp^{\cell}a , \susp^{\cell}b).\]
We observe that those $(n+1)$\-cells are parallel again by
Property~\ref{susp-2}.

\paragraph{Fullness.}
To finish the induction, it remains to show that for every Batanin tree $B$ and
full $n$\-sphere $A = (a,b)$ of $\free_n \Pos{B}$, the sphere $\susp^{\type} A$
is full in $\free_{n+1}\Pos{\susp B}$. To show that, we first let
\begin{align*}
  a &= \cell_n\free_n(s_n^B)(a_0) & b &= \cell_n\free_n(t_n^B)(b_0),
\end{align*}
where the support of $a_0$ and $b_0$ contains all positions of $\partial_n B$.
Then $\susp^{\type}A$ consists of the cells
\begin{align*}
  a' &= \susp^{\cell}(\cell_n\free_n(s_n^B)(a_0)) \\
     &= \cell_n\susp\free_n(s_n^B)(\susp^{\cell}(a_0)) \\
     &= \cell_n\free_n(\susp s_n^B)(\susp^{\cell}(a_0)) \\
     &= \cell_n\free_n(s_{n+1}^{\susp B})(\susp^{\cell}(a_0)) \\
  b' &= \cell_n\free_n(t_{n+1}^{\susp B})(\susp^{\cell}(b_0)),
\end{align*}
so it remains to show that when the support of $a$ contains all positions of a
tree $B$, then the support of $\susp^{\cell}a$ contains all positions of
$\susp B$. More generally, it suffices to prove that for every computad $C$ and
every cell $c\in \cell_n C$,
\[\supp(\susp^{\cell}(c)) = \supp(c) \cup \set{v_-,v_+}.\]
This statement can be easily shown by structural induction on cells.

\paragraph{Infinite-dimensional computads.}
This concludes the induction on $n\in \N$. Compatibility of the suspension
functors with the forgetful functors allows us to define a functor
\[
  \susp : \comp \to \compbip
\]
sending a computad $C = (C_n)_{n\in \N}$ to the computad with components
\begin{align*}
  (\susp C)_0 &= \{v_-,v_+\} &
                               (\susp C)_{n+1} &= \susp C_n,
\end{align*}
and with basepoints $\var v_-$ and $\var v_+$. By construction, this functor
commutes with the suspension operation on globular sets, in that the following
square commutes:
\[\begin{tikzcd}
  \comp & \compbip \\
  \glob & \bipointed
          \arrow["\susp", from=2-1, to=2-2]
          \arrow["\susp", from=1-1, to=1-2]
          \arrow["\free", from=2-1, to=1-1]
          \arrow["\freebip"', from=2-2, to=1-2]
\end{tikzcd}\]
Moreover, Property~\ref{susp-2} shows that the natural transformations
$\susp^{\cell}$ can be combined to a natural transformation
\[ \susp^{\cell} : \susp \cell \Rightarrow \cellbip\susp. \]

\begin{example}\label{ex:susp-composite}
  The suspension of the tree \(B_{n,k,m}\) of Example~\ref{ex:ps-comp} is the
  tree \(\susp B_{n,k,m} = B_{n+1,k+1,m+1}\). Moreover, applying the natural
  transformation \(\susp^{\cell}\) to the cell \(\unbcomp_{n,k,m}\) defined in
  Section~\ref{sec:cell-compositions} yields the cell
  \[
    \susp^{\cell}(\unbcomp_{n,k,m}) = \unbcomp_{n+1,k+1,m+1}.
  \]
  It follows that any cell of the form \(\unbcomp_{n,k,m}\) can be obtained by
  iteratively suspending one of the form \(\unbcomp_{r,0,s}\).
\end{example}

\subsection{Hom \texorpdfstring{$\omega$}{ω}-categories}
We can define bipointed \(\omega\)\-categories similarly to bipointed globular
sets and bipointed computads. Alternatively, the adjunction
\[
  \freebip : \bipointed \rightleftarrows \compbip : \cellbip
\]
gives rise to a monad \(\Tbip\) on bipointed globular sets, sending a bipointed
globular set $(X, x_-,x_+)$ to the globular set $TX$ with basepoints $\var x_-$
and $\var x_+$. It can be shown that its category of algebras is the category
$\catbip$ of bipointed $\omega$\-categories.

Whiskering the natural transformation \(\susp^{\cell}\) on the right with
$\free$ and using compatibility of the suspension with the functor \(\free\), we
get a natural transformation
\[
  \susp^T : \susp T \Rightarrow \Tbip\susp.
\]
Equivalently, by the mate
correspondence~\cite[Proposition~2.1]{kelly_review_1974}, we get a natural
transformation
\[
  \deloop^T = (\deloop\Tbip\kappa)\circ (\deloop \susp^T \deloop) : T\deloop
  \Rightarrow \deloop \Tbip.
\]
We will show that this \((\deloop,\deloop^{T})\) is a morphism of monads,
meaning that the following diagrams commute:
\[\begin{tikzcd}
	T\deloop && \deloop\Tbip \\
	& \deloop
	\arrow["{\deloop^T}", Rightarrow, from=1-1, to=1-3]
	\arrow["\eta\deloop", Rightarrow, from=2-2, to=1-1]
	\arrow["\deloop\eta"', Rightarrow, from=2-2, to=1-3]
\end{tikzcd}\qquad
\begin{tikzcd}
	TT\deloop & {T\deloop \Tbip} & \deloop\Tbip\Tbip \\
	T\deloop && \deloop\Tbip
	\arrow["{T\deloop^T}", Rightarrow, from=1-1, to=1-2]
	\arrow["{\deloop^T\Tbip}", Rightarrow, from=1-2, to=1-3]
	\arrow["\deloop\mu", Rightarrow, from=1-3, to=2-3]
	\arrow["\mu\deloop"', Rightarrow, from=1-1, to=2-1]
	\arrow["{\deloop^T}"', Rightarrow, from=2-1, to=2-3]
\end{tikzcd}\]
By the mate correspondence, commutativity of those diagrams is equivalent to the
commutativity of the following ones:
\[\begin{tikzcd}
	{\susp T } && \Tbip\susp \\
	& \susp
	\arrow["\susp\eta", Rightarrow, from=2-2, to=1-1]
	\arrow["\eta\susp"', Rightarrow, from=2-2, to=1-3]
	\arrow["{\susp^T}", Rightarrow, from=1-1, to=1-3]
\end{tikzcd}\qquad
\begin{tikzcd}
	{\susp TT} & {\Tbip \susp T} & \Tbip\Tbip\susp \\
	{\susp T} && \Tbip\susp
	\arrow["{\susp^TT}", Rightarrow, from=1-1, to=1-2]
	\arrow["{\Tbip \susp^T}", Rightarrow, from=1-2, to=1-3]
	\arrow["\mu\susp", Rightarrow, from=1-3, to=2-3]
	\arrow["{\susp \mu}"', Rightarrow, from=1-1, to=2-1]
	\arrow["{\susp^T}"', Rightarrow, from=2-1, to=2-3]
\end{tikzcd}\]
The left one commutes, since $\susp^{\cell}$ preserves generators. The right one
is obtained from the following diagram by whiskering on the right with $\free$:
\begin{equation}\label{cd:sigma-cell-monad-law}
  \begin{tikzcd}
	{\susp T\cell} & {\Tbip \susp \cell} & \Tbip\cellbip\susp \\
	{\susp \cell} && \cellbip\susp
	\arrow["{\susp^{\cell}\free\cell}"{outer sep = 5pt}, Rightarrow, from=1-1, to=1-2]
	\arrow["{\Tbip \susp^{\cell}}", Rightarrow, from=1-2, to=1-3]
	\arrow["\cellbip\varepsilon\susp", Rightarrow, from=1-3, to=2-3]
	\arrow["{\susp \cell\varepsilon}"', Rightarrow, from=1-1, to=2-1]
	\arrow["{\susp^{\cell}}"', Rightarrow, from=2-1, to=2-3]
  \arrow["\cellbip \susp \varepsilon"', Rightarrow, from=1-2, to=2-3]
\end{tikzcd}
\end{equation}
This diagram commutes by naturality of \(\Sigma^{\cell}\) and commutativity of
the following diagram
\[\begin{tikzcd}[column sep = large]
  \susp \free \cell \dar[equals]{}
  \rar[Rightarrow]{\susp\varepsilon} &
  \susp \\
  \freebip \susp \cell \rar[Rightarrow]{\freebip \susp^{\cell}} &
  \freebip \cellbip \susp \uar[Rightarrow]{\varepsilon \susp}
\end{tikzcd}\]
This can be checked easily by showing that both sides agree on generators.

\begin{definition}
  Consider a bipointed $\omega$\-category
  \((X, {\alpha : TX \to X}, x_-, x_+)\). Its \emph{hom $\omega$\-category} is
  the $\omega$\-category that consists of the globular set $\deloop (X,x_,x_+)$
  and the structure morphism
  \[
    T\deloop (X,x_,x_+) \xrightarrow{\deloop^T}
    \deloop \Tbip X \xrightarrow{\deloop \alpha}
    \deloop (X,x_,x_+)
  \]
\end{definition}
\noindent
As shown by Street~\cite{street_formal_1972}, and explained by
Leinster~\cite[Theorem~6.1.1]{leinster_higher_2004},
this definition extends to a functor
\[\deloop : \catbip \to \omega\Cat,\]
since \((\deloop,\deloop^{T})\) is a morphism of monads. By construction, this
functor makes the following solid square commute:
\[
  \begin{tikzcd}
    \catbip & \omega\Cat \\
    \bipointed & \glob
    \arrow["\deloop", from=2-1, to=2-2]
    \arrow["\deloop", from=1-1, to=1-2]
    \arrow["{U^T}", from=1-2, to=2-2]
    \arrow["{U^{\Tbip}}"', from=1-1, to=2-1]
    \arrow["\susp", from=2-2, to=2-1, dashed, bend left]
    \arrow["\susp"', from=1-2, to=1-1, dashed, bend right]
  \end{tikzcd}
\]

\begin{example}
  In light of the \(\omega\)\-category structure on the hom globular set, the
  last equation of Example~\ref{ex:susp-composite} can be understood in
  particular as saying that the vertical composite of two \(2\)\-cells \(x,y\)
  in an \(\omega\)\-category \(\wcat X\) can be described as the composite of
  the \(1\)\-cells \(x,y\) in the hom \(\omega\)\-category
  \(\deloop(\wcat X,\src x,\tgt x)\). Similar reasoning also shows that the
  associators and unitors for the composition of \(1\)\-cells in the hom
  \(\omega\)\-category \(\deloop(\wcat X, \src x,\tgt x)\) yields through the
  natural transformation \(\susp^{\cell}\) associators and unitors for the
  vertical composition in \(\wcat X\).
\end{example}

\paragraph{Suspension.}
The category of bipointed \(\omega\)\-categories is cocomplete being the
category of algebras of a finitary monad on a locally presentable category, and
the vertical functors in this commuting square are monadic. Therefore, by the
adjoint lifting theorem~\cite{johnstone_adjoint_1975}, the hom
\(\omega\)\-category functor also admits a left adjoint \(\susp\) that commutes
with the free functors \(F^{\Tbip}\) and \(F^T\). We will show that restricting
the left adjoint to the subcategory of computads gives the suspension functor
defined before, in the sense that the following diagram commutes up to
isomorphism:
\[
  \begin{tikzcd}
    \catbip & \omega\Cat \\
    \compbip & \comp
    \arrow["\susp", from=2-1, to=2-2]
    \arrow["\susp", from=1-1, to=1-2]
    \arrow["K^{\Tbip}", from=2-1, to=1-1]
    \arrow["{K^{T}}"', from=2-2, to=1-2]
  \end{tikzcd}
\]
Since the suspension of an \(\omega\)\-category is defined as a left adjoint,
commutativity of this square amounts to showing for every computad \(C\) and
every bipointed \(\omega\)\-category \((Y,y_-,y_+)\) that there exists a natural
bijection
\[
  \psi : \catbip(K^{\Tbip}(\susp C), Y) \xrightarrow{\sim}
  \omega\Cat(K^T C, \Omega Y)
\]
We claim that such a bijection is given by
\[
  \psi(f) = \Omega (f \circ \susp_C^{\cell}) : K^T C \to \Omega Y
\]
First of all, we need to check that \(\psi(f)\) is a morphism of
\(\omega\)\-categories for every morphism of bipointed \(\omega\)\-categories
\(f\). This amounts to commutativity of the exterior of the following diagram
\[\begin{tikzcd}
	{T\cell C} & {T \deloop \susp \cell C} & {T\deloop \cellbip \susp C} & {T \deloop Y} \\
	&& {\deloop \Tbip \cellbip \susp C} & {\Omega \Tbip Y} \\
	{\cell C} & {\deloop \susp \cell C} & {\deloop \cellbip \susp C} & {\deloop Y}
	\arrow[equals, from=1-1, to=1-2]
	\arrow[equals, from=3-1, to=3-2]
	\arrow["{T\deloop \susp^{\cell}}", from=1-2, to=1-3]
	\arrow["{T\deloop f}", from=1-3, to=1-4]
	\arrow["{\deloop^T}", from=1-4, to=2-4]
	\arrow[from=2-4, to=3-4]
	\arrow["{\deloop f}"', from=3-3, to=3-4]
	\arrow["{\deloop \susp^{\cell}}"', from=3-2, to=3-3]
	\arrow["{\cell(\varepsilon_C)}"', from=1-1, to=3-1]
	\arrow["{\deloop \Tbip f}", from=2-3, to=2-4]
	\arrow["{\deloop \cellbip(\varepsilon_{\susp C})}"', from=2-3, to=3-3]
	\arrow["{\deloop^T}"', from=1-3, to=2-3]
\end{tikzcd}\]
where the unnamed morphism \(\Tbip Y\to Y\) is the \(\Tbip\)\-algebra structure of
\(Y\). The left square in this diagram commutes, being the transpose of diagram
\eqref{cd:sigma-cell-monad-law}. The top square commutes by naturality, and the
bottom square commutes by \(f\) being an algebra morphism. Therefore,
\(\psi(f)\) is well\-defined.

To show that \(\psi\) is bijective, we use that morphisms out of a computad are
determined by their values on
generators~\cite[Corollary~6.5]{dean_computads_2024}. Assume that
\(\psi(f) = \psi(g)\) for a pair of morphisms of bipointed
\(\omega\)\-categories. Then \(f\) and \(g\) agree on the \(0\)\-dimensional
generators being bipointed, and they agree on positive-dimensional ones by
\(\Omega\susp^{\cell}\) being a bijection on the generators.
Therefore, \(f = g\). For surjectivity, given a morphism
\(h : K^T C\to \Omega Y\), we define inductively a sequence of morphisms
\(f_n : K_n^T (\susp C)_n \to Y\) by letting \(f_0\) being the unique morphism
determined by \(f(v_\pm) = y_\pm\) and \(f_{n+1}\) being the morphism induced
by \(f_n\) and the function
\[
  V_{n+1}^{\susp C} = V_n^C \xrightarrow{\var} \cell_n C
  \xrightarrow{h} (\Omega Y)_n \hookrightarrow Y_{n+1}
\]
Finally, we let \(f \colon K^T(\susp C) \to Y\) be the morphisms out of the
colimits induced by the \(f_n\). Then by construction \(\psi(f)\) and \(h\)
agree on generators, so they must be equal.

\subsection{The hom of a computad.}
In this section, we show that the hom \(\omega\)\-category functor preserves the
property of being free on a computad. This is contrary to the case of strict
\(\omega\)\-categories: consider for instance the strict \(\omega\)\-category
freely generated on the computad \(C_{\eh}\) with only one object
\(x\), and two generators \(a,b \colon \id x \to \id x\) of dimension \(2\). By
the Eckmann-Hilton argument, the composition of \(1\)\-cells in
\(\deloop(C_{\eh},x,x)\) is commutative, which prevents it from being
free on a computad.

We will show that \(\deloop K^{\Tbip} C\) is free on a computad when \(C\) is a
bipointed computad for weak \(\omega\)\-categories. To do so, we construct a
computad \(\deloop C\) on which it is free. The generators of \(\deloop C\) are
precisely the \emph{indecomposable cells} of the \(\omega\)\-category
\(\deloop K^T C\):
\begin{definition}
  An \(n\)\-cell of an \(\omega\)\-category \(\wcat X = (X, \alpha : TX \to X)\)
  is \emph{indecomposable} when it is not of the form
  \(\alpha(\coh{B}{A}{\free\sigma})\) for any Batanin tree \(B\), full
  \((n-1)\)\-type \(A\) and morphism \(\sigma : \Pos B\to X\). We will denote
  the set of indecomposable \(n\)\-cells by \(X_n^{\ind}\).
\end{definition}

\noindent Using the indecomposable cells of an \(\omega\)\-category \(\wcat X\),
we can build a computad \(C_{\wcat X}\) together with a morphism
\(\sigma_{\wcat X} : K^T C_{\wcat X}\to {\wcat X}\). We start inductively by
letting \(C_{\wcat X, -1}\) the unique \(-1\)\-computad and
\(\sigma_{\wcat X,-1} : K_{-1}^TC_{\wcat X, -1} \to \wcat X\) the unique
morphism from the initial \(\omega\)\-category. We then define the
\(n\)\-computad \( C_{\wcat X, n} \) to be the triple
\((C_{\wcat X, n-1}, V^{\wcat X}_n , \phi_n^{\wcat X})\) defined via the
following pullback square
\begin{equation}
  \begin{tikzcd}[column sep = huge]\label{cd:cnx-def}
	{V_n^{\wcat X}} & {X^{\ind}_n} \\
	& {X_n } \\
	{\type_{n-1}(C_{\wcat X,n-1})} & {\Par_{n-1}(X)}
	\arrow["{(\sigma_{\wcat X,n-1},\sigma_{\wcat X,n-1})}", from=3-1, to=3-2]
	\arrow[hook, from=1-2, to=2-2]
	\arrow["{(\src,\tgt)}", from=2-2, to=3-2]
	\arrow["{\sigma^{\ind}_{\wcat X, n,V}}", dashed, from=1-1, to=1-2]
	\arrow["{\phi_n^{\wcat X}}"', dashed,  from=1-1, to=3-1]
  \arrow[from=1-1, to=3-2, phantom, "\ulcorner"{rotate=180,very near start}]
  \end{tikzcd}
\end{equation}
where \(\Par_{n-1}X\) is the set of pairs of parallel \((n-1)\)\-cells of \(X\).
The morphism \(\sigma_{\wcat X,n} : K_n^T C_{\wcat X, n}\to {\wcat X}\) is
the one determined by \(\sigma_{\wcat X, n-1}\) and the function
\(V_n^{\wcat X}\to X_n\) in the diagram above. Finally, we define
\(C_{\wcat X}\) to consist of the computads \(C_{\wcat X,n}\), and we let
\(\sigma_{\wcat X}\) to be the morphism out of the colimit induced by the
morphisms \(\sigma_{\wcat X,n}\).

\begin{proposition}
  For every bipointed computad \(C\), the \(\omega\)\-category
  \(\Omega K^{\Tbip}C\) is free on a computad.
\end{proposition}
\begin{proof}
  Let \(\wcat X = \Omega K^{\Tbip}C\). We will show that the morphism
  \(\sigma_{\wcat X}\) is an isomorphism, or equivalently that it is bijective
  on all cells. We proceed by strong induction on the dimension of cells. Let
  therefore \(n \ge 0\) and suppose that the result holds for all \(k < n\).
  Since the inclusion \(K^T_{n-1}C_{n-1}\to K^TC\) is the identity on cells of
  dimension at most \((n-1)\), the bottom morphism of the
  square~\eqref{cd:cnx-def} must be an isomorphism, and hence the top one should
  be as well. This implies that \(\sigma_{\wcat X}\) is injective on generator
  \(n\)\-cells.

  This shows that \(\sigma_{\wcat X}\) is injective on \(0\)\-cells. If
  \(n > 0\), then we can see that \(\sigma_{\wcat X}\) sends a coherence
  \(n\)\-cell \(c = \coh{B}{A}{\tau}\) to the decomposable \(n\)\-cell:
  \[\begin{split}
    \sigma_{\wcat X}(c)
      &= \sigma_{\wcat X}
        (\coh{B}{A}{\varepsilon \circ \free \tau^{\dagger}}) \\
      &= (\sigma_{\wcat X} \circ
        \cell(\varepsilon))(\coh{B}{A}{\free \tau^{\dagger}}) \\
      &= (\alpha_{\wcat X}\circ \cell\free \sigma_{\wcat X})
        (\coh{B}{A}{\free \tau^{\dagger}}) \\
      &= \alpha_{\wcat X}(\coh{B}{A}{\free (\sigma_{\wcat X}\tau^{\dagger})})
  \end{split}\]
where \(\tau^\dagger : \Pos B\to \cell\free C\) the transpose of \(\tau\)
under the \(\free \dashv \cell\) adjunction. As \(\sigma_{\wcat X}\) sends
generators to indecomposable cells, it can not send a generator and a
coherence to the same cell.

  Finally, to show that \(\sigma_{\wcat X}\) is
  injective on \(n\)\-cells, it remains to show that it is injective on
  coherence cells. For that, let \(c = \coh{B}{A}{\tau}\) and
  \(c' = \coh{B'}{A'}{\tau'}\) and suppose that
  \[
    \sigma_{\wcat X}(c) = \sigma_{\wcat X}(c').
  \]
  From the definition of \(\alpha_{\wcat X}\), we can compute further that
  \[\begin{split}
    \sigma_{\wcat X}(c)
      &= (\deloop\cellbip(\varepsilon_C)(\deloop\Tbip \kappa_{\cellbip C})
         (\deloop \susp^T_{\deloop \cellbip C}))
        (\coh{B}{A}{\free (\sigma_{\wcat X}\tau^{\dagger})}) \\
      &= (\deloop\cellbip(\varepsilon_C)(\deloop\Tbip \kappa_{\cellbip C}))
        (\coh{\susp B}{\susp^{\type}A}
        {\susp\free(\sigma_{\wcat X}\tau^{\dagger})}) \\
      &= \coh{\susp B}{\susp^{\type}A}{\varepsilon_C
        \free\kappa\circ\susp\free(\sigma_{\wcat X}\tau^{\dagger})} \\
      &= \coh{\susp B}{\susp^{\type}A}{\varepsilon_C\circ
        \free(\kappa \circ \susp (\sigma_{\wcat X}\tau^{\dagger}))}
  \end{split}\]
so from the equality above and injectivity of the constructor \(\coherence\),
we get that \(\susp B = \susp B'\), that \(\susp^{\type}A = \susp^{\type}A'\), and
that the corresponding morphisms agree. The suspension is injective
on trees by injectivity of \(\bt\), so \(B = B'\). The equality of
morphisms implies that
\(\sigma_{\wcat X}\tau^{\dagger} = \sigma_{\wcat X}\tau^{\prime\dagger}\) by
transposing along the adjunction \(\free\vdash\cell\) and
\(\susp\vdash\deloop\). By structural induction on the \(n\)\-cells, we may
assume that \(\sigma_{\wcat X}\) is injective on the \(n\)\-cells in the
image of \(\tau^\dagger\) from which we conclude that
\(\tau^{\dagger} = \tau^{\prime\dagger}\) and hence that \(\tau = \tau'\).
Finally by structural induction, we can see that \(\susp^{\cell}\) and
\(\susp^{\type}\) are monic, so \(A = A'\) and \(c = c'\).

  To show that \(\sigma_{\wcat X}\) is surjective on \(n\)\-cells, we proceed by
  structural induction on the cells. Since the function
  \(\sigma_{\wcat X,n,V}^{\ind}\) is bijective, every indecomposable
  \(n\)\-cell is in the image of \(\sigma_{\wcat X}\). Suppose therefore that
  \(c \in X_n\) can be decomposed in the form
  \[
    c = \alpha_{\wcat X}(\coh{B}{A}{\free \tau}) =
    \coh{\susp B}{\susp^{\type}A}
    {\varepsilon_C \circ \free (\kappa \circ \susp \tau)}.
  \]
  By structural induction, we may assume that for every position
  \(p\in \Pos[k]{B}\), there exists a cell \(\tau'(p)\in \cell_n(C_{\wcat X})\)
  such that
  \[\tau(p) = \sigma_{\wcat X}(\tau'(p))\]
  By injectivity of \(\sigma\) on cells of dimension at most \(k\), those cells
  are unique and they can be assembled into a morphism
  \(\tau' : \Pos B\to \cell(C_{\wcat X})\) such that
  \(\tau = \sigma_{\wcat X} \tau'\). Using this morphism and the calculations
  above, we see that
  \[
    c =
    \alpha_{\wcat X}(\coh{B}{A}{\free (\sigma_{\wcat X} \tau')}) =
    \sigma_{\wcat X}(\coh{B}{A}{\varepsilon\circ \free \tau'}),
  \]
  so \(\sigma_{\wcat X}\) is surjective on \(n\)\-cells.
\end{proof}

\begin{corollary}
  The assignment \(C \mapsto C_{\deloop K^{\Tbip} C}\) extends to a functor
  \[
    \deloop : \compbip \to \comp
  \]
  making the following two diagrams commute up to isomorphism:
  \begin{align*}
    \begin{tikzcd}[ampersand replacement=\&]
      \catbip \& \omega\Cat \\
      \compbip \& \comp
      \arrow["{K^{\Tbip}}", from=2-1, to=1-1]
      \arrow["{K^T}"', from=2-2, to=1-2]
      \arrow["\deloop", from=1-1, to=1-2]
      \arrow["\deloop"', from=2-1, to=2-2]
    \end{tikzcd}
    &&
       \begin{tikzcd}[ampersand replacement=\&]
         \compbip \& \comp \\
         \bipointed \& \glob
         \arrow["\deloop", from=1-1, to=1-2]
         \arrow["\cellbip"', from=1-1, to=2-1]
         \arrow["\cell", from=1-2, to=2-2]
         \arrow["\deloop"', from=2-1, to=2-2]
       \end{tikzcd}
  \end{align*}
\end{corollary}
\begin{proof}
  Since the comparison functor \(K^{T}\) is fully faithful, for every morphism
  of computads \(\tau : C \to D\), we define \(\deloop \tau\) to be the unique
  morphism making the following diagram commute:
  \[
    \begin{tikzcd}
      K^{T}\deloop C
      \ar[r,"\sigma","\sim"']
      \ar[d,dashed,"K^{T}\deloop \tau"']
      & \deloop K^{\Tbip}C
      \ar[d,"\deloop K^{\Tbip}\tau"]\\
      K^{T}\deloop D
      \ar[r,"\sigma","\sim"']
      & \deloop K^{\Tbip}D
    \end{tikzcd}
  \]
\end{proof}

\begin{remark}
  In general, the hom computad functor does not preserve finiteness. For
  example, consider the computad \(\free \disk 0\) with only one
  \(0\)\-dimensional generator \(x\). The hom computad
  \(\deloop(\free\disk 0,x,x)\) has a countable number of \(0\)\-dimensional
  generators, corresponding to the identity of \(x\) and possible ways it can be
  composed.
\end{remark}


%% file: 4-opposites.tex

\section{Opposites}
\label{sec:opp}
An important feature of ordinary category theory is the duality stemming from
the existence of opposite categories. This feature extends to higher categories,
where we may define opposites by reversing the direction of all cells in certain
dimensions. In this section, we will define the opposite of a globular set, a
computad, and an $\omega$\-category with respect to a set of dimensions
$w \subseteq \N_{>0}$, using the same technique as in the previous section. We
will then show that the formation of opposites in all those cases gives rise to
an action of the Boolean group
\[ G = \mathcal{P}(\N_{>0}) \cong \Z_2^{\N_{>0}}\]
of subsets of the positive integers with respect to symmetric difference. This
group is clearly isomorphic to the group of functions $\N_{>0} \to \Z_2$ with
pointwise multiplication, where each subset is identified with its indicator
function. Abusing notation we will identify a subset $w$ with its indicator
function, and write $w(n)$ for the value of the indicator function at
$n\in\N_{>0}$.

\subsection{The opposite of a globular set}
The group $G$ acts on the category $\G$ of globes by swapping the source and
target inclusions. More precisely, an element $w\in G$ acts as the
identity-on-objects functor
\[\op_w \colon \G \to \G \]
given on the generating morphisms by
\begin{align*}
  \op_w(s_n) &=  \begin{cases}
    t_n & \text{if } n+1 \in w, \\
    s_n & \text{if } n+1 \not\in w,
  \end{cases} &
  \op_w(t_n) &=  \begin{cases}
    s_n & \text{if } n+1 \in w, \\
    t_n & \text{if } n+1 \not\in w.
  \end{cases}
\end{align*}
The functor $\op_{\emptyset}$ is clearly the identity functor. Moreover, for
every pair of elements $w, w' \in G$, we can easily check that
\[\op_w\op_{w'} = \op_{ww'}\]
so the assignment $w \mapsto \op_w$ is a group homomorphism $G \to \Aut(\G)$.
Since the group $G$ is Abelian, this action extends to an action on the category
$\glob$ of globular sets by precomposition
\begin{align*}
  \op \colon G &\to \Aut(\glob) \\
  \op_w(X) &= X \circ \op_w.
\end{align*}
The opposite $\op_wX$ of a globular set $X$ therefore has the same cells as $X$,
with the source and target of $n$\-cells reversed for $n\in w$.

Since pasting diagrams are bipointed by their $0$\-source and $0$\-target
inclusions, it will be useful to further extend this action to an action on
bipointed globular sets
\[
  \op \colon G\to \Aut(\bipointed)
\]
by letting $\op_w$ take a bipointed globular set $(X,x_-,x_+)$ to the opposite
globular set $\op_w X$ with the same basepoints when $1\not\in w$, and with the
basepoints swapped otherwise.

\begin{lemma}\label{lem:op-susp}
  For every $w\in G$, there exists a natural isomorphism
  \[
    \op_w^{\susp} \colon
    \susp \circ\ \op_{w-1} \Rightarrow
    \op_w \circ \susp
  \]
  where $w-1\in G$ is the sequence defined by $(w-1)(n) = w(n+1)$. Moreover,
  $\op^{\susp}_\emptyset$ is the identity natural transformation, and for
  every pair of elements $w, w'\in G$, the following diagram commutes:
  \[\begin{tikzcd}[column sep = large, row sep = small]
      {\susp\op_{w-1}\op_{w'-1}} &
      {\op_w\susp \op_{w'-1}} &
      {\op_{w}\op_{w'}\susp} \\
      {\susp\op_{ww'-1}} && {\op_{ww'}\susp}
      \arrow[Rightarrow, "{\op_w^{\susp}\op_{w'-1}}", from=1-1, to=1-2]
      \arrow[Rightarrow, "{\op_w\op^{\susp}_{w'}}", from=1-2, to=1-3]
      \arrow[Rightarrow, "{\op_{ww'}^{\susp}}"', from=2-1, to=2-3]
      \arrow[equals, from=1-1, to=2-1]
      \arrow[equals, from=1-3, to=2-3]
  \end{tikzcd}\]
\end{lemma}

\begin{proof}
  For every globular set $X$, the bipointed globular sets $\susp \op_{w-1}X$
  and $\op_w \susp X$ have the same sets of cells. Moreover, the source and
  target of an $n$\-cell in both of them agree when $n>2$: they are given by
  the target and source functions of $X$ respectively when $n\in w$, and they
  are given by the source and target functions of $X$ when $n\not\in w$. The
  source and target of a $1$\-cell in the first one are given by $v^-$ and
  $v^+$ respectively, while in the latter it is given by those when
  $1\not\in w$, and by $v^+$ and $v^-$ when $1\in w$. Therefore, we may define
  an isomorphism of globular sets
  \[
    \op^{\susp}_{w,X} \colon \susp \op_{w-1}X \to \op_w \susp X
  \]
  to be the identity on positive-dimensional cells, and to be given on
  $0$\-cells by
  \[
    \op_{w,X}^{\susp}(v_\pm) =
    \begin{cases}
      \op_{w,X}^{\susp}(v_\mp), &\text{if } 1\in w, \\
      \op_{w,X}^{\susp}(v_\pm), &\text{if } 1\not\in w.
    \end{cases}
  \]
  Since $\op_w$ reverses the basepoints if and only if $1\in w$, we see that
  this is a morphism of bipointed globular sets. Naturality of these morphisms
  follows easily by the fact that it is the identity of positive-dimensional
  cells. Finally, the claimed diagram commutes for $w, w'\in G$: both
  morphisms are identity on positive-dimensional cells, they are the identity
  on $0$\-cells when $1\in w\cap w'$ or $1\not\in w\cup w'$, and they swap the
  two $0$\-cells otherwise.
\end{proof}

\begin{lemma}\label{lem:op-wedge}
  For every $w\in G$ and $n\in\N$, there exists a natural isomorphism
  \[
    \op_w^{\vee} \colon
    \bigvee_{i=1}^n \circ\ \swap_{w(1)} \circ (\op_w)^n
    \Rightarrow \op_w \circ \bigvee_{i=1}^n
  \]
  where $\swap_0$ is the identity of $(\bipointed)^n$, while $\swap_1$ is the
  automorphism
  \[\swap_1(X_1,\dots,X_n) = (X_n,\dots,X_1).\]
  Moreover, $\op_\emptyset^\vee$ is the identity natural transformation, and
  for every pair of elements $w, w'\in G$, the following diagram
  commutes:
  \[\begin{tikzcd}
      {\bigvee \circ (\op_{ww'})^n\circ\swap_{ww'(1)}} &
      {\op_{ww'}\circ \bigvee} \\
      {\bigvee \circ (\op_{w})^n\circ (\op_{w'})^n \circ
          \swap_{w(1)}\circ\swap_{w'(1)}} \\
      {\bigvee \circ (\op_{w})^n\circ\swap_{w(1)}\circ
          (\op_{w'})^n \circ \swap_{w'(1)}} \\
      {\op_w\circ\bigvee\circ (\op_{w'})^n \circ \swap_{w'(1)}} &
      {\op_w\op_w' \circ \bigvee}
      \arrow[equals, from=1-1, to=2-1]
      \arrow[equals, from=2-1, to=3-1]
      \arrow["{\op_{w}^{\vee}}", from=3-1, to=4-1]
      \arrow["{\op_w(\op_{w'}^{\vee})}", from=4-1, to=4-2]
      \arrow["{\op_{ww'}^{\vee}}", from=1-1, to=1-2]
      \arrow[equals, from=1-2, to=4-2]
  \end{tikzcd}\]
\end{lemma}

\begin{proof}
  Fix $n\in \N$ and $w\in G$ and let $X_1,\dots,X_n$ be bipointed globular
  sets and suppose first that $1\not\in w$, so that the basepoints of $X_i$
  and $\op_w X_i$ agree. The functor $\op_w$ on globular sets preserves
  $\disk 0$, and it preserves colimits, being an equivalence of categories.
  Therefore, there exists a natural isomorphism of globular sets
  \[
      \op_w^\vee \colon \bigvee_{i=1}^n (\op_w X_i) \to
      \op_w \left(\bigvee_{i=1}^n X_i\right),
  \]
  that can be easily seen to preserve the basepoints. Moreover, since $\op_w$
  preserves the cells of a globular set, and colimits of globular sets are
  computed pointwise, we may take $\op_w^\vee$ to be the identity.

  Suppose now that $1\in w$, so that the functor $\op_w$ swaps the basepoints.
  Using that $\op_w$ preserves colimits and $\disk 0$, we see that
  $\op_w(\bigvee_{i=1}^n X_i)$ is the colimit of the following diagram.
  \[\begin{tikzcd}[column sep = large]
      {\op_w X_1} && \dots && {\op_w X_n} \\
      & {\disk 0} && {\disk 0}
      \arrow["{x_{n,-}}"', from=2-4, to=1-5]
      \arrow["{x_{1,+}}", from=2-2, to=1-1]
      \arrow["{x_{2,-}}"', from=2-2, to=1-3]
      \arrow["{x_{n-1,+}}", from=2-4, to=1-3]
  \end{tikzcd}\]
  On the other hand, $\bigvee_{i=n}^1 \op_w X_i$ is the colimit of the
  following diagram:
  \[\begin{tikzcd}[column sep = large]
      {\op_w X_n} && \dots && {\op_w X_1} \\
      & {\disk 0} && {\disk 0}
      \arrow["{x_{1,+}}"', from=2-4, to=1-5]
      \arrow["{x_{n,-}}", from=2-2, to=1-1]
      \arrow["{x_{n-1,+}}"', from=2-2, to=1-3]
      \arrow["{x_{2,-}}", from=2-4, to=1-3]
  \end{tikzcd}\]
  By symmetry of pushouts, we get a natural isomorphism of globular sets
  \[
    \op_w^\vee \colon \bigvee_{i=n}^1 (\op_w X_i) \to
    \op_w \left(\bigvee_{i=1}^n X_i\right)
  \]
  that can be easily seen to preserve the basepoints. Since colimits are
  computed object-wise, this isomorphism is given level-wise by the symmetry of
  pushouts.

  Knowing how those isomorphisms are defined pointwise, we can easily deduce
  that the claimed diagram commutes for every pair $w, w' \in G$. If
  $1\not\in w\cup w'$, then both sides of the diagram are identities. If
  $1\in w\cap w'$ again both are identities, since the symmetry of the
  pushout squares to the identity. Finally, when $1\in ww'$, then both sides
  are given by the symmetry of the pushout, so they agree.
\end{proof}

Using those lemmas, we can deduce that pasting diagrams are closed under the
formation of opposites: we define recursively on the Batanin tree $B$ for every
$w\in G$ the \emph{$w$\-opposite Batanin tree} $\op_wB$ by the formula
\[
  \op_w(\bt{B_1,\dots,B_n})=
  \btlist(\swap_{w(1)}[\op_{w-1}B_1,\dots,\op_{w-1}B_n]),
\]
where $\swap_0$ is the identity of the set of lists, while $\swap_1$ reverses a
list
\[ \swap_1[B_1,\dots,B_n] = [B_n,\dots,B_1].\]
The opposite tree realizes the opposite pasting diagram, in the sense that there
exists an isomorphism of bipointed globular sets
\[ \op_w^B \colon \Pos{\op_wB} \to \op_w(\Pos B). \]
We can define this isomorphism recursively on $B = \bt{B_1,\dots,B_n}$ to be
the following composite
\begin{align*}
  \Pos{\op_wB}
    &= \bigvee_{i=1}^n\swap_{w(1)}(\susp \Pos{\op_{w-1}B_i}) \\
    &\xrightarrow{\bigvee\swap_{w(1)}\susp \op_{w-1}^{B_i}}
      \bigvee_{i=1}^n\swap_{w(1)}(\susp \op_{w-1}\Pos{B_i}) \\
    &\xrightarrow{\bigvee\swap_{w(1)}\op^{\susp}_{w}}
      \bigvee_{i=1}^n\swap_{w(1)}(\op_w\susp\Pos{B_i}) \\
    &\xrightarrow{\op_w^{\vee}}
      \op_w\left(\bigvee_{i=1}^n \susp\Pos{B_i}\right) = \op_w(\Pos{B}).
\end{align*}

\begin{lemma}\label{lem:op-tree}
  The isomorphism $\op_\emptyset^B$ is the identity for every tree $B$, and
  for any $w, w'\in G$, the following diagram of isomorphisms commutes:
  \[\begin{tikzcd}[column sep = large]
      {\Pos{\op_w\op_{w'}B}} & {\op_w(\Pos{\op_{w'}B})}
      & {\op_w\op_{w'}\Pos{B}} \\
      {\Pos{\op_{ww'}B}} && {\op_{ww'}\Pos{B}}
      \arrow["{\op_{w}^{\op_{w'}B}}", from=1-1, to=1-2]
      \arrow["{\op_w(\op_{w'}^{B})}", from=1-2, to=1-3]
      \arrow[equals, from=1-3, to=2-3]
      \arrow[equals, from=1-1, to=2-1]
      \arrow["{\op_{ww'}^{B}}"', from=2-1, to=2-3]
  \end{tikzcd}\]
\end{lemma}
\begin{proof}
  This lemma is an easy induction on $B$, using naturality of the isomorphisms
  in Lemmas~\ref{lem:op-susp} and~\ref{lem:op-wedge}, and of the commuting
  diagrams there.
\end{proof}

\begin{lemma}\label{lem:op-bdry}
  For every $w\in G$, $k\in \N$ and Batanin tree $B$,
  \[\op_w\partial_k = \partial_k \op_w. \]
  Moreover, the following equations hold
  \[
    \begin{array}{r@{ = }l|r@{ = }l}
      \multicolumn{2}{c|}{k+1\in w}
      & \multicolumn{2}{c}{k+1\notin w}\\[1.2em]
      \hline \multicolumn{2}{c|}{} & \multicolumn{2}{|c}{} \\
      \op_w(t_k^B) \circ \op_w^{\partial_k B}
      &\op_w^B \circ\ s_k^{\op_w B}
      &\op_w(s_k^B) \circ \op_w^{\partial_k B}
      &\op_w^B \circ\ s_k^{\op_w B} \\[1.2em]
      \op_w(s_k^B) \circ \op_w^{\partial_k B}
      & \op_w^B \circ\ t_k^{\op_w B}
      & \op_w(t_k^B) \circ \op_w^{\partial_k B}
      & \op_w^B \circ\ t_k^{\op_w B}
    \end{array}
  \]
\end{lemma}
\begin{proof}
  We proceed by induction on $k$. For $k = 0$ both $\op_w \partial_k B$ and
  $\partial_k \op_w B$ are equal to the disk $D_0$, and the equations state that
  $\op_w^B$ preserves the basepoints. Suppose therefore that the result is true
  for some $k\in \N$ to prove that it also holds for $k+1$. Letting
  $B = \bt{B_1,\dots,B_n}$, we see that
  \begin{align*}
    \op_w\partial_{k+1} B
    &= \op_w(\bt{\partial_kB_1,\dots,\partial_{k}B_n}) \\
    &= \btlist(\swap_{w(1)}[\op_{w-1}\partial_kB_1,\dots,\op_{w-1}\partial_kB_n]) \\
    &= \btlist(\swap_{w(1)} [\partial_k \op_{w-1} B_1,\dots,\partial_k \op_{w-1}B_n]) \\
    &= \partial_{k+1}\btlist(\swap_{w(1)} [\op_{w-1} B_1,\dots,\op_{w-1}B_n]) \\
    &= \partial_{k+1} \op_w B
\end{align*}
  by the inductive hypothesis.

  We will prove the first equation in the case that $k+1 \in w$ and $1\in w$.
  The other equation and the rest of the cases follow by the same argument.
  By the inductive hypothesis, we may assume that for $1\le i \le n$, the
  following square commutes
  \[\begin{tikzcd}
      {\Pos{\op_{w-1} B_i}} & {\op_{w-1}(\Pos{B_i})} \\
      {\Pos{\partial_{k-1}\op_{w-1} B_i}} \\
      {\Pos{\op_{w-1}\partial_{k-1} B_i}} &
      {\op_{w-1}\Pos{\partial_{k-1} B_i}}
      \arrow["{s_{k-1}^{\op_{w-1}B_i}}", from=2-1, to=1-1]
      \arrow["{\op_{w-1}^{B_i}}", from=1-1, to=1-2]
      \arrow[equals, from=2-1, to=3-1]
      \arrow["{\op_{w-1}^{\partial_{k-1} B_i}}"', from=3-1, to=3-2]
      \arrow["{\op_{w-1}(t_{k-1}^{B_i})}"', from=3-2, to=1-2]
  \end{tikzcd}\]
  Applying the suspension functor and then the wedge sum from $n$ to $1$, we
  get that the left square below commutes. Naturality of the isomorphisms in
  Lemmas~\ref{lem:op-susp} and~\ref{lem:op-wedge} then imply that the right
  square below also commutes.
  \[\begin{tikzcd}[column sep = tiny]
      {\bigvee\susp\Pos{\op_{w-1} B_i}}
      & {\bigvee\susp\op_{w-1}\Pos{B_i}}
      & {\op_w\bigvee\susp\Pos{B_i}} \\
      {\bigvee\susp\Pos{\partial_{k-1} \op_{w-1} B_i}} \\
      {\bigvee\susp\Pos{\op_{w-1}\partial_{k-1} B_i}}
      & {\bigvee\susp\op_{w-1}\Pos{\partial_{k-1}B_i}}
      & {\op_{w}\bigvee\susp\Pos{\partial_{k-1}B_i}}
      \arrow[equals, from=2-1, to=3-1]
      \arrow["{\bigvee\susp s_{k-1}^{\op_{w-1}B_i}}"', from=2-1, to=1-1]
      \arrow[from=1-1, to=1-2]
      \arrow[from=3-1, to=3-2]
      \arrow["{\bigvee\susp\op_{w-1}(t_{k-1}^{B_i})}"{description},
          from=3-2, to=1-2]
      \arrow[from=3-2, to=3-3]
      \arrow[from=1-2, to=1-3]
      \arrow["{\op_w\bigvee\susp t_{k-1}^{B_i}}"{description},
          from=3-3, to=1-3]
  \end{tikzcd}\]
  The outer part of the diagram though is precisely the square:
  \[\begin{tikzcd}
      {\Pos{\op_{w} B}} & {\op_{w}(\Pos{B})} \\
      {\Pos{\partial_{k}\op_{w} B}} \\
      {\Pos{\op_{w}\partial_{k} B}} &
      {\op_{w}\Pos{\partial_{k} B}}
      \arrow["{s_{k}^{\op_{w}B}}", from=2-1, to=1-1]
      \arrow["{\op_{w}^{B}}", from=1-1, to=1-2]
      \arrow[equals, from=2-1, to=3-1]
      \arrow["{\op_{w}^{\partial_{k} B}}"', from=3-1, to=3-2]
      \arrow["{\op_{w}(t_{k}^{B})}"', from=3-2, to=1-2]
    \end{tikzcd}\]
  whose commutativity amounts to the first equation.
\end{proof}

\begin{remark}
  Batanin trees famillially represent the free strict \(\omega\)\-category monad
  \(T^{\str}\) on the category of globular sets, in the sense that
  \[
    (T^{\str}X)_n = \coprod_{\dim B \le n} \glob(\Pos{B}, X)
  \]
  as shown by Leinster~\cite{leinster_higher_2004}. Compatibility of the
  opposites with Batanin trees allows us to define an invertible morphism of
  monads \(\op_w^{\str}\colon T^{\str}\op_w \Rightarrow \op_w T^{\str}\) by the
  formula
  \[
    \op_{w,X}^{\str}(B, f \colon \Pos B\to \op_w X) =
    (\op_w B, \op_w f\circ \op_w^B)
  \]
  hence providing an alternative definition of the opposites of a strict
  \(\omega\)\-category.
\end{remark}

\begin{example}
  The disk trees \(D_{n}\) defined in Example~\ref{ex:ps-disk} are invariant
  under the action of opposites. In other words for every \(w\in G\), we have
  \(\op_{w} D_{n} = D_{n}\). However, in general, the isomorphism
  \(\op_w^{D_n}\) is not the identity. The family of trees \(B_{n,k,m}\) of
  Example~\ref{ex:ps-comp} satisfies the following equality with respect to the
  opposites:
  \[
    \op_{w}(B_{n,k,m}) =
    \begin{cases}
      B_{m,k,n} & \text{if \(k+1 \in w\)} \\
      B_{n,k,m} & \text{otherwise.}
    \end{cases}
  \]
\end{example}
\begin{example}
  As a richer example, consider the tree \(B = \bt{\bt{\bt{\bt{}},\bt{}},\bt{}}\).
  We illustrate the various opposites of this tree on the corresponding globular
  pasting diagram with the following figure:
  \[
    \begin{tikzpicture}
      \coordinate (A) at (0,2);
      \coordinate (B) at (4,4);
      \coordinate (C) at (4,0);
      \coordinate (D) at (8,2);

    \node (a) at (A) {
      \(\begin{tikzcd}[ampersand replacement=\&]
        \bullet
        \ar[r, bend left=50]
        \ar[r, bend left=25, phantom, "\Downarrow\Rrightarrow\Downarrow"]
        \ar[r]
        \ar[r, bend right=25, phantom, "\Downarrow"]
        \ar[r, bend right=50]
        \& \bullet
        \ar[r]
        \& \bullet
      \end{tikzcd}\)
    };

    \node (b) at (B) {
      \(\begin{tikzcd}[ampersand replacement=\&]
        \bullet
        \ar[r, bend left=50]
        \ar[r, bend left=25, phantom, "\Downarrow"]
        \ar[r]
        \ar[r, bend right=25, phantom, "\Downarrow\Rrightarrow\Downarrow"]
        \ar[r, bend right=50]
        \& \bullet
        \ar[r]
        \& \bullet
      \end{tikzcd}\)
    };

    \node (c) at (C) {
      \(\begin{tikzcd}[ampersand replacement=\&]
        \bullet
        \ar[r]
        \& \bullet
        \ar[r, bend left=50]
        \ar[r, bend left=25, phantom, "\Downarrow\Rrightarrow\Downarrow"]
        \ar[r]
        \ar[r, bend right=25, phantom, "\Downarrow"]
        \ar[r, bend right=50]
        \& \bullet
      \end{tikzcd}\)
    };

    \node (d) at (D) {
      \(\begin{tikzcd}[ampersand replacement=\&]
        \bullet
        \ar[r]
        \& \bullet
        \ar[r, bend left=50]
        \ar[r, bend left=25, phantom, "\Downarrow"]
        \ar[r]
        \ar[r, bend right=25, phantom, "\Downarrow\Rrightarrow\Downarrow"]
        \ar[r, bend right=50]
        \& \bullet
      \end{tikzcd}\)
    };
    \draw[|->] (a) -- (b) node[midway, above] {\(\op_2\)};
    \draw[|->] (a) -- (c) node[midway, above] {\(\op_1\)};
    \draw[|->] (a) -- (d) node[midway, above] {\(\op_{1,2}\)};
    \draw[|->] (b) -- (d) node[midway, above] {\(\op_1\)};
    \draw[|->] (c) -- (d) node[midway, above] {\(\op_2\)};
  \end{tikzpicture}
  \]
\end{example}

\subsection{The opposite of a computad}
The opposite of a computad is defined similarly to the opposite of a globular
set by swapping the source and target of its generators. To define this action,
we fix an element $w\in G$, that will be omitted from the notation, and
define recursively on the dimension $n \in \N$,
an endofunctor and two natural transformation
\begin{align*}
  \op &\colon \comp_n \to \comp_n \\
  \op^{\cell} &\colon \cell_n \Rightarrow\cell_n \op  \\
  \op^{\type} &\colon \type_n \Rightarrow \type_n \op
\end{align*}
satisfying the following properties:
\begin{enumerate}[label={\textsc{(op{\footnotesize\arabic*})}}]
\item\label{op-1} forming opposites commutes with the forgetful functors, and
  the inclusion of globular sets into computads
  \[\begin{tikzcd}
    {\comp_{n+1}} & {\comp_{n+1}} \\
    {\comp_n} & {\comp_n}
                \arrow["{\op}", from=1-1, to=1-2]
                \arrow["{u_{n+1}}"', from=1-1, to=2-1]
                \arrow["{u_{n+1}}", from=1-2, to=2-2]
                \arrow["{\op}"', from=2-1, to=2-2]
  \end{tikzcd}\qquad\qquad
  \begin{tikzcd}
    \glob & \glob \\
    {\comp_n} & {\comp_n}
                \arrow["{\free_n}"', from=1-1, to=2-1]
                \arrow["{\op}"', from=2-1, to=2-2]
                \arrow["{\op}", from=1-1, to=1-2]
                \arrow["{\free_n}", from=1-2, to=2-2]
  \end{tikzcd}\]
\item\label{op-2} the natural transformations are compatible with the boundary
  natural transformation
  \[\begin{tikzcd}
    {\cell_{n+1}}
    && {\cell_{n+1}\op} \\
    {\type_n u_{n+1}}
    & {\type_n \op u_{n+1}}
    &
      {\type_n u_{n+1}\op}
      \arrow["{\ty_{n+1}}"', Rightarrow, from=1-1, to=2-1]
      \arrow["{\op^{\type}u_{n+1}}"', Rightarrow, from=2-1, to=2-2]
      \arrow["{\op^{\cell}}", Rightarrow, from=1-1, to=1-3]
      \arrow["{\ty_{n+1}\op}", Rightarrow, from=1-3, to=2-3]
      \arrow[equals, from=2-2, to=2-3]
  \end{tikzcd}\]
\item\label{op-3} the natural transformation $\op^{\type}$ swaps the two cells
  of a sphere when $n+1\in w$ and leaves them unchanged otherwise, in the sense
  that the following diagrams commute for $i = 1, 2$
  \[\begin{tikzcd}
    {\type_n} & {\type_n\op} \\
    {\cell_n} & {\cell_n\op}
                \arrow["{\pr_i}"', Rightarrow, from=1-1, to=2-1]
                \arrow[""{name=0, anchor=center, inner sep=0},
                "{\op^{\cell}}"', Rightarrow, from=2-1, to=2-2]
                \arrow[""{name=1, anchor=center, inner sep=0},
                "{\op^{\type}}",Rightarrow, from=1-1, to=1-2]
                \arrow["{\pr_i}\op", Rightarrow, from=1-2, to=2-2]
                \arrow["{n+1\not\in w}"{description}, draw=none, from=1, to=0]
  \end{tikzcd}\qquad\qquad
  \begin{tikzcd}
    {\type_n} & {\type_n\op} \\
    {\cell_n} & {\cell_n\op} \arrow["{\pr_i}"', Rightarrow, from=1-1,
                to=2-1] \arrow[""{name=0, anchor=center, inner sep=0},
                "{\op^{\cell}}"', Rightarrow, from=2-1, to=2-2]
                \arrow[""{name=1, anchor=center, inner sep=0},
                "{\op^{\type}}", Rightarrow, from=1-1, to=1-2]
                \arrow["{\pr_{2-i}\op}", Rightarrow, from=1-2, to=2-2]
                \arrow["{n+1\in w}"{description}, draw=none, from=1, to=0]
  \end{tikzcd}\]
\item\label{op-4} the natural transformation $\op^{\cell}$ preserves
  generators, in that for every globular set $X$ and $x\in X_n$, we have that
  \[
    \op_{\free X}^{\cell}(\var x) = \var x.
  \]
\item\label{op-5} the natural transformation $\op^{\type}$ preserves fullness,
  in that for every full $n$\-sphere $A$ of $\free\Pos{B}$, the $n$\-sphere
  \[A' = \type_{n}\free_{n}(\op^B)^{-1}(\op^{\type}(A))\] of
  $\free\Pos{\op B}$ is also full.
\end{enumerate}
As a base case, we define $\op $ and $\op^{\type}$ to be the identities of
$\comp_{-1}$ and $\type_{-1}$ respectively. Let therefore $n\in \N$ and suppose
inductively that data as above has been defined for all natural numbers less
than $n$, satisfying the given properties.

\paragraph{Computads.}
First we will define the action of \(\op\) on all \(n\)-computads. Let
$C = (C_{n-1}, V_n^C, \phi_n^C)$ be an $n$\-computad. The opposite computad
$\op C$ consists of the opposite computad $\op C_{n-1}$, the same set of
generators $V_n^C$, and the attaching function
\[
  \phi_n^{\op C} \colon V_n^C \xrightarrow{\phi_n^C} \type_{n-1} C_{n-1}
  \xrightarrow{\op^{\type}} \type_{n-1}(\op C_{n-1}).
\]
By Properties~\ref{op-3} and~\ref{op-4}, we can easily deduce that $\op$
commutes with the inclusion $\free_n$ on objects, while it clearly commutes with
the forgetful functors $u_{n}$ by definition.

\paragraph{Cells and morphisms.}
We will then define $\op$ on morphisms \(\sigma\) of \(n\)-computads of target
\(C\), together with the component of the natural transformation $\op^{\cell}$
at \(C\) mutually recursively. For a generator $v\in V_n^C$, we let
\[\op^{\cell}(\var v) = \var v,\]
and we observe that
\[
  \ty_{n}\op^{\cell}(\var v) = \op^{\type}\ty_{n}(\var v).
\]
Given a coherence cell $c = \coh{B}{A}{\tau}$ of $C$, we may assume that
recursively that $\op(\tau)$ has been defined, and let
\begin{align*}
  A' &= \type_{n-1}\free_{n-1}(\op^B)^{-1}(\op^{\type}(A)) \\
  \op^{\cell}(c) &= \coh{\op B}{A'}{\op(\tau)\circ \free_n(\op^B)}
\end{align*}
We then observe again that the boundary of this cell is given by
\begin{align*}
  \ty_{n}\op^{\cell}(c)
  &= \type_{n-1}(\op\tau_{n-1})(\op^{\type} A) \\
  &= \op^{\type}(\type_{n-1}(\tau_{n-1})(A)) \\
  &= \op^{\type} \ty_n(c).
\end{align*}
Finally, for a morphism $\sigma = (\sigma_{n-1}, \sigma_V) : D\to C$, we define
assume that $\op^{\cell}$ has been defined on cells of the form $\sigma_V(v)$
for $v\in V_n^D$ and define
\[
  \op(\sigma) = (\op\sigma_{n-1}, \op^{\cell}\circ\ \sigma_V) : \op D\to
  \op C
\]
This is a well-defined morphism of computads by the observation on the boundary
of the cells $\op^{\cell}(c)$, i.e. by Properties~\ref{op-2}.

It follows immediately from the definition that $\op$ commutes with the
forgetful functor $u_{n}$ on morphisms as well. Using that $\op^{\cell}$
preserves generators, we can also deduce that $\op$ commutes with the
inclusion $\free_n$ on morphisms as well. Therefore, we have shown
Properties~\ref{op-1},~\ref{op-2}~and~\ref{op-4} so far.

\paragraph{Naturality.}
We will now show that $\op$ is a functor and that $\op^{\cell}$ is natural.
For that, we fix a morphism of $n$\-computads $\sigma \colon C\to D$, and we
proceed recursively to show that the following square commutes
\[\begin{tikzcd}[column sep = large]
  {\cell_n C} & {\cell_nD} \\
  {\cell_n\op C} & {\cell_n\op D}
                    \arrow["{\cell_n\sigma}", from=1-1, to=1-2]
                    \arrow["{\op^{\cell}}"', from=1-1, to=2-1]
                    \arrow["{\op^{\cell}}", from=1-2, to=2-2]
                    \arrow["{\cell_n\op\sigma}"', from=2-1, to=2-2]
\end{tikzcd}\]
and that for all morphism $\tau \colon E\to C$,
\[\op\sigma \circ \op\tau = \op(\sigma\circ \tau).\]
By definition of $\op\sigma$, the square above commutes when restricted to
generators. Moreover, for a coherence cell $c = \coh{B}{A}{\tau}$, we see that
\begin{align*}
  \op^{\cell}\circ \cell_n(\sigma)(c)
  &= \op^{\cell}(\coh{B}{A}{\sigma\circ\tau}) \\
  &= \coh{B}{A'}{\op(\sigma\circ\tau) \circ \free(\op^B)} \\
  &= \coh{B}{A'}{\op(\sigma)\circ \op(\tau)\circ\free(\op^B)} \\
  &= \cell_n(\op\sigma)(\coh{B}{A'}{\op(\tau)\circ\free(\op^B)}) \\
  &= \cell_n(\op\sigma)\circ\op^{\cell}(\coh{B}{A}{\tau})
\end{align*}
where $A'$ is defined as above. Given arbitrary $\tau \colon E\to C$, we may
assume that the square commutes when restricted to the image of $\tau_{V}$. By
the inductive hypothesis, \(\op\) preserves composition of morphisms of
\((n-1)\)-computads. Hence it suffices to show the equality above for the
generators of E. We recall the definition of the composition of morphisms of
\(n\)\-computads given in~\cite[Section~3.1]{dean_computads_2024}:
\[
  (\sigma_{n-1},\sigma_{v}) \circ (\tau_{n-1},\tau_{V}) = (\sigma_{n-1}\circ
  \tau_{n-1}, \cell_{n}(\sigma)\circ \tau_{v})
\]
Using this definition, we have:
\begin{align*}
  (\op(\sigma\circ\tau))_V
  &= \op^{\cell}\circ \cell_n(\sigma)\circ\tau_V  \\
  &= \cell_n(\op\sigma)\circ \op_n^{\cell}\circ \tau_V \\
  &= (\op(\sigma)\circ \op(\tau))_V.
\end{align*}
Therefore, $\op$ is a functor and $\op^{\cell}$ is natural.

\paragraph{Spheres.}
The natural transformation $\op^{\type}$ is completely determined by
Property~\ref{op-3}. Indeed, for an $n$\-computad $C$ and for a sphere
$(a,b) \in \type_nC$, we are forced to define
\[
  \op^{\type}(a,b) =
  \begin{cases}
    (\op^{\cell}b, \op^{\cell}a),
    &\text{if } n+1\in w \\
    (\op^{\cell}a, \op^{\cell}b),
    &\text{if } n+1\not\in w \\
  \end{cases}
\]
Property~\ref{op-2} shows us that those $\op^{\cell}a$ and $\op^{\cell}b$
have the same source and target, so that this assignment is well-defined. It is
clearly natural by naturality of $\op^{\cell}$.

\paragraph{Fullness.}
To finish the recursive definition, it remains to show that for every Batanin
tree $B$ and every $n$\-sphere $A = (a,b)$ of $\free_n\Pos{B}$, the $n$\-sphere
\[
  A' = \type_{n}\free_{n}(\op^B)^{-1}(\op^{\type}(A))
\]
of $\free_n\Pos{\op B}$ is also full. We will show that in the case that
$n+1\in w$, the other case being similar. By assumption, we may write
\begin{align*}
  a &= \cell_n\free_n(s_n^B)(a_0) & b &= \cell_n\free_n(t_n^B)(b_0).
\end{align*}
For $n$\-cells $a_0,b_0$ of $\free_n\Pos{\partial_n B}$ whose support contains
all positions of $\partial_n B$. Then we have that $A' = (a',b')$ where
\begin{align*}
  a' &= \cell_{n}\free_{n}((\op^B)^{-1} \circ t_n^B)(\op^{\cell}b_0) \\
  b' &= \cell_{n}\free_{n}((\op^B)^{-1} \circ t_n^B)(\op^{\cell}a_0)
\end{align*}
By Lemma~\ref{lem:op-bdry}, we can rewrite those cells as
\begin{align*}
  a'&= \cell_{n}\free_{n}(s_n^{\op B})
      (\cell_n\free_n(\op^{\partial_nB})^{-1}(\op^{\cell}b_0)) \\
  b'&= \cell_{n}\free_{n}(t_n^{\op B})
      (\cell_n\free_n(\op^{\partial_nB})^{-1}(\op^{\cell}a_0)).
\end{align*}
Using the definition of the support and that $\op^{\cell}$ preserves
generators, we may show recursively that
\[\supp(\op^{\cell}(c)) = \supp(c) \]
for every cell $c$. Moreover, isomorphisms of computads induce bijections on the
support of cells, so the support of the cells
\begin{align*}
  \cell_n\free_n(\op^{\partial_nB})^{-1}(\op^{\cell}b_0) \\
  \cell_n\free_n(\op^{\partial_nB})^{-1}(\op^{\cell}a_0)
\end{align*}
must contain all positions of $\partial_n\op B$. Therefore, $A'$ is full.

\begin{lemma}\label{lem:op-fin-comp}
  For every $n\in \N$, the endofunctor $\op_\emptyset$ on $n$\-computads is the
  identity, and so are the natural transformations $\op_\emptyset^{\cell}$ and
  $\op_\emptyset^{\type}$. Moreover, for any pair of elements $w, w'\in G$,
  \[\op_w \op_{w'} = \op_{ww'}\]
  and the following diagrams commute.
  \[\begin{tikzcd}[column sep = huge]
    {\cell_n} & {\cell_n\op_{w'}} & {\cell_n\op_w\op_{w'}} \\
    {\cell_n} && {\cell_n\op_{ww'}}
                 \arrow["{\op_{ww'}^{\cell}}", Rightarrow, from=2-1, to=2-3]
                 \arrow[equals, from=1-1, to=2-1]
                 \arrow[equals, from=1-3, to=2-3]
                 \arrow["{\op_{w'}^{\cell}}", Rightarrow, from=1-1, to=1-2]
                 \arrow["{\op_{w}^{\cell}\op_{w'}}", Rightarrow, from=1-2, to=1-3]
  \end{tikzcd}\]
\[\begin{tikzcd}[column sep = huge]
  {\type_n} & {\type_n\op_{w'}} & {\type_n\op_w\op_{w'}} \\
  {\type_n} && {\type_n\op_{ww'}}
               \arrow["{\op_{ww'}^{\type}}", Rightarrow, from=2-1, to=2-3]
               \arrow[equals, from=1-1, to=2-1]
               \arrow[equals, from=1-3, to=2-3]
               \arrow["{\op_{w'}^{\type}}", Rightarrow, from=1-1, to=1-2]
               \arrow["{\op_{w}^{\type}\op_{w'}}", Rightarrow, from=1-2, to=1-3]
\end{tikzcd}\]
In particular, $\op_w$, $\op_w^{\cell}$ and $\op_w^{\type}$ are invertible
with inverses $\op_w$, $\op_w^{\cell}\op_w$ and $\op_w^{\type}\op_w$
respectively.
\end{lemma}
\begin{proof}
  We proceed inductively on $n\in\N$, since the result holds trivially for
  $n = -1$. Since $\op_\emptyset$ and $\op_\emptyset^{\type}$ are identities for
  $(n-1)$\-computads, we see that
  \[\op_\emptyset C = C\]
  for every $n$\-comptutad $C$. Using Lemma~\ref{lem:op-bdry}, we can then show
  mutually recursively for an $n$\-computad $C$ that
  \begin{align*}
    \op_\emptyset \sigma &= \sigma &
                                     \op_\emptyset^{\cell} c &= c
  \end{align*}
  for every morphism $\sigma \colon D\to C$ and every $n$\-cell $c$ of $C$.
  Using then that $\op_\emptyset^{\type}$ is defined using
  $\op_\emptyset^{\cell}$, we see that $\op_\emptyset^{\type}$ must be the
  identity as well.

  Let now $w, w'\in G$ and $C = (C_{n-1}, V_n^C, \phi_n^C)$ an $n$\-computad.
  Then the \(n\)\-computad $\op_w\op_{w'}C$ consists of the $(n-1)$\-computad
  \[
    \op_w\op_{w'}C_{n-1} = \op_{ww'} C_{n-1},
  \]
  the same set of generators, and the attaching function
  \[
    \phi_n^{\op_w\op_{w'}C} = \op_{w,\op_{w'}C}^{\type}\circ
    \op_{w',C}^{\type}\circ\ \phi_n^C = \op_{ww',C}^{\type}\circ\ \phi_n^C =
    \phi_n^{\op_{ww'}C}.
  \]
  Hence, $\op_w\op_{w'}$ and $\op_{ww'}$ agree on $n$\-computads. Fixing a
  computad $C$, we can show that they also agree on morphisms with target $C$
  mutually inductively to recursively to showing that the claimed diagram for
  $\op_w^{\cell}$ commutes. The commutative diagram from $\op_w^{\type}$ then
  follows from the one for $\op_w^{\cell}$.
\end{proof}

Having defined the opposite of an $n$\-computad for every $n\in\N$, in a way
that is compatible with the forgetful functors $u_n$, we get for every
\(w \in G\), a functor
\[ \op \colon \comp \to \comp \] sending a computad $C = (C_n)_{n\in\N}$ to
the computad
\[\op C = (\op C_n)_{n\in\N}\]
and acts similarly on morphisms. Property~\ref{op-1} shows that $\op$ is
compatible with the inclusion functors $\free$ in that
\[\begin{tikzcd}
  \comp & \comp \\
  \glob & \glob
          \arrow["{\op}", from=1-1, to=1-2]
          \arrow["{\op}", from=2-1, to=2-2]
          \arrow["\free"', from=2-2, to=1-2]
          \arrow["\free", from=2-1, to=1-1]
\end{tikzcd}\]
commutes. Moreover, combining Properties~\ref{op-2} and~\ref{op-3}, we see that
the natural transformations $\op^{\cell}$ give rise to a natural
transformation
\[ \op^{\cell} \colon \op \cell \Rightarrow \cell \op. \] The following
lemma is an easy consequence of Lemma~\ref{lem:op-fin-comp}.

\begin{lemma}\label{lem:op-comp}
  The functor $\op_\emptyset \colon \comp\to \comp$ is the identity functor, and
  the natural transformations $\op_\emptyset^{\cell}$ is the identity of
  $\cell$. Moreover, for any pair $w, w' \in G$,
  \[ \op_w \op_{w'} = \op_{ww'}\] and the following diagrams commute.
  \[\begin{tikzcd}
    {\op_w\op_{w'}\cell} & {\op_w\cell\op_{w'}} & {\cell \op_w\op_{w'}} \\
    {\op_{ww'}\cell} && {\cell\op_{ww'}}
                        \arrow[equals, from=1-1, to=2-1]
                        \arrow[equals, from=1-3, to=2-3]
                        \arrow[Rightarrow, "{\op_{ww'}^{\cell}}"', from=2-1, to=2-3]
                        \arrow[Rightarrow, "{\op_w\op_{w'}^{\cell}}", from=1-1, to=1-2]
                        \arrow[Rightarrow, "{\op_w^{\cell}\op_{w'}}", from=1-2, to=1-3].
  \end{tikzcd}\]
In particular, each $\op_w$ is invertible with inverse itself, and
$\op_w^{\cell}$ is invertible with inverse $\op_w\op_w^{\cell}\op_w$.
\end{lemma}

\begin{example}\label{ex:op-whisk}
  The family of cells \(\unbcomp_{n,k,m} \in T\Pos{B_{n,k,m}}\) defined in
  Section~\ref{sec:cell-compositions} satisfies the following identities under
  the action of composites:
  \[
    \op^{\cell}_{w}(\unbcomp_{n,k,m}) =
    \begin{cases}
      \unbcomp_{m,k,n} & \text{if \(k+1 \in w\)} \\
      \unbcomp_{n,k,m} & \text{otherwise.}
    \end{cases}
  \]
  We defined \(\unbcomp_{n,k,m}\) with \(m > n\) by analogy, informally relying
  on the reader's ability to consrtuct it from the case where \(n < m\). The
  construction of opposites give us a way to make this argument formal.
\end{example}

\subsection{The opposite of an \texorpdfstring{$\omega$}{ω}-category.}
So far, we have defined the opposite of a globular set, a pasting diagram and a
computad. To extend those definitions and define the opposite of an
$\omega$\-category, we consider the mate
\begin{align*}
  \op^T &: T\op_w \Rightarrow \op T \\
  \op^T &= (\op^{\cell} \free)^{-1}
\end{align*}
of the natural transformation \(\op^{\cell}\free\) under the adjunction
\(\op \dashv \op\). We will show that the functor $\op : \glob\to \glob$
together with the natural transformation \(\op^{T}\) is a morphism of monads
from $T$ to $T$. This amounts to commutativity of the following two diagrams
\[\begin{tikzcd}
	{T\op} && {\op T} \\
	& {\op}
	\arrow[Rightarrow, "{\eta\op}", from=2-2, to=1-1]
	\arrow[Rightarrow, "{\op\eta}"', from=2-2, to=1-3]
	\arrow[Rightarrow, "{\op^T}", from=1-1, to=1-3]
\end{tikzcd}\qquad
\begin{tikzcd}
	{TT\op} & {T\op T} & {\op TT} \\
	{T\op} && {\op T}
	\arrow[Rightarrow, "{\op^T}", from=2-1, to=2-3]
	\arrow[Rightarrow, "{\mu\op}"', from=1-1, to=2-1]
	\arrow[Rightarrow, "{\op\mu}", from=1-3, to=2-3]
	\arrow[Rightarrow, "{T\op^T}", from=1-1, to=1-2]
	\arrow[Rightarrow, "{\op^TT}", from=1-2, to=1-3]
\end{tikzcd}\]
The left one is the assertion that $\op^{\cell}$ preserves generators, which
we have already shown. The right one is obtained from the following diagram by
whiskering on the right with $\free$, and then replacing $\op^{\cell}\free$ with
its inverse.
\[\begin{tikzcd}
	{\op T\cell} && {\cell\op\free\cell} & {T\op\cell} & {T\cell\op} \\
	{\op \cell} &&&& {\cell\op}
	\arrow[Rightarrow, "{\op^{\cell}\free\cell}"{outer sep=4pt},
    from=1-1, to=1-3]
	\arrow[equals, from=1-3, to=1-4]
	\arrow[Rightarrow, "{T\op^{\cell}}", from=1-4, to=1-5]
	\arrow[Rightarrow, "{\cell\varepsilon\op}", from=1-5, to=2-5]
	\arrow[Rightarrow, "{\op\cell\varepsilon}"', from=1-1, to=2-1]
	\arrow[Rightarrow, "{\op^{\cell}}"', from=2-1, to=2-5]
  \arrow[Rightarrow, "{\cell\op\varepsilon}"', from = 1-3, to = 2-5]
\end{tikzcd}\]
The left square commutes by naturality of $\op^{\cell}$, while the one on the
right is obtained from the following square by whiskering with \(\cell\) on he
left.
\[\begin{tikzcd}[column sep = large]
  \op \free \cell \dar[equals]{} \rar[Rightarrow]{\op \varepsilon} &
  \op \\
  \free \op \cell \rar[Rightarrow]{\free \op^{\cell}} &
  \free \cell \op \uar[Rightarrow]{\varepsilon \op}
\end{tikzcd}\]
Finally we can check that this square commutes by showing that two sides agree
on every generator.

\begin{definition}
  The opposite of an $\omega$\-category $(X,\alpha : TX\to X)$ with respect to
  some $w\in G$ is the $\omega$\-category consisting of the globular set
  $\op_w X$ and the structure morphism
  \[
    T\op_w X \xrightarrow{\op_{w,X}^T} \op_w TX \xrightarrow{\op_w\alpha}\op_wX
  \]
\end{definition}

\noindent The construction of the opposite of an $\omega$\-category is
well-defined and gives rise to endofunctors
\[\op_w : \omega\Cat \to \omega\Cat\]
for every $w\in G$ as shown by Street~\cite{street_formal_1972}, and explained
by Leinster~\cite[Theorem~6.1.1]{leinster_higher_2004}. Moreover, the following
lemma - an immediate consequence of Lemma~\ref{lem:op-comp} - shows that those
endofunctors are invertible and give rise to an action
\[ \op : G\to \Aut(\omega\Cat) \]
of $G$ on the category of $\omega$\-categories.

\begin{lemma}
  The natural transformation $\op_\emptyset^T$ is the identity, and for any
  ${w, w'\in G}$ the following diagram commutes:
  \[\begin{tikzcd}[column sep = large]
    {T\op_w\op_{w'}} & {\op_wT\op_{w'}} & {\op_{w}\op_{w'}T} \\
    {T\op_{ww'}} && {\op_{ww'}T}
    \arrow[Rightarrow, "{\op_w^T\op_{w'}}", from=1-1, to=1-2]
    \arrow[Rightarrow, "{\op_{w}\op_{w'}^T}", from=1-2, to=1-3]
    \arrow[equals, from=1-3, to=2-3]
    \arrow[equals, from=1-1, to=2-1]
    \arrow[Rightarrow, "{\op_{ww'}^T}"', from=2-1, to=2-3]
  \end{tikzcd}\]
\end{lemma}
The commutative diagram showing that $\op^T_w$ is a morphism of monads implies
in particular that the components of the natural transformation $\op_w^{\cell}$
are morphisms of free $\omega$\-categories, so it can also be seen as a natural
isomorphism
\[  \op_w^K : \op_w K^T \Rightarrow K^T\op_w \]
The opposite functors on globular sets, computads and $\omega$\-categories are
therefore related by the following five squares
\[\begin{tikzcd}
	\omega\Cat & \omega\Cat & \omega\Cat & \omega\Cat & \comp & \comp \\
	\glob & \glob & \comp & \comp & \glob & \glob
	\arrow["{U^T}"', from=1-1, to=2-1]
	\arrow["{U^T}", from=1-2, to=2-2]
	\arrow["{\op_w}", from=1-1, to=1-2]
	\arrow["{\op_w}", from=2-1, to=2-2]
	\arrow["{K^T}", from=2-3, to=1-3]
	\arrow["{K^T}"', from=2-4, to=1-4]
	\arrow["{\op_w}", from=2-3, to=2-4]
	\arrow["{\op_w}", from=1-3, to=1-4]
	\arrow["\op_w^K", shorten <=4pt, shorten >=4pt, Rightarrow, from=2-3, to=1-4]
	\arrow["\free"', from=2-6, to=1-6]
	\arrow["\free", from=2-5, to=1-5]
	\arrow["{\op_w}", from=2-5, to=2-6]
	\arrow["{\op_w}", from=1-5, to=1-6]
\end{tikzcd}\]
\[\begin{tikzcd}
	\comp & \comp & \omega\Cat & \omega\Cat \\
	\glob & \glob & \glob & \glob
	\arrow["{\op_w}", from=2-1, to=2-2]
	\arrow["{\op_w}", from=1-1, to=1-2]
	\arrow["\cell"', from=1-1, to=2-1]
	\arrow["\cell", from=1-2, to=2-2]
	\arrow["{\op^{\cell}_w}", shorten <=4pt, shorten >=4pt,
    Rightarrow, from=2-1, to=1-2]
	\arrow["{\op_w}", from=2-3, to=2-4]
	\arrow["{\op_w}", from=1-3, to=1-4]
	\arrow["{F^T}", from=2-3, to=1-3]
	\arrow["{F^T}"', from=2-4, to=1-4]
	\arrow["{\op_w^F}", shorten <=4pt, shorten >=4pt,
    Rightarrow, from=2-3, to=1-4]
\end{tikzcd}\]
where
\begin{align*}
  \op_w^F &= \op_w^K\free & \op_w^{\cell} &= U^T\op_w^K.
\end{align*}
In particular, the opposite of an $\omega$\-category that is free on a globular
set or computad is again free on the opposite of the underlying globular set or
the opposite computad, up to natural isomorphism.

\begin{example}
  In the light of the action of opposites on \(\omega\)\-categories,
  specialising Example~\ref{ex:op-whisk} to the cell \(\unbcomp_{1,0,2}\) can be
  interpreted as stating that the left whiskering of a cell \(y\) by a cell
  \(x\) in an \(\omega\)\-category \(\wcat X\) can be described as the right
  whiskering of \(y\) by \(x\) in the opposite \(\omega\)\-category
  \(\op_{1}(\wcat X)\).
\end{example}

\begin{remark}
  Since the free \(\omega\)\-category monad \(T\) is isomorphic to Leinster's
  initial contractible globular operad, it is equipped with a natural
  transformation \(\alpha \colon T \Rightarrow T^{\str}\) to the free strict
  \(\omega\)\-category monad admitting a contraction. An alternative way to
  construct the morphism of monads \(\op^T \colon T\op \Rightarrow \op T\),
  would be to show that the composite
  \[
    \op T \op
    \xRightarrow{\op \alpha \op} \op T^{\str} \op
    \xRightarrow{\op \op^{\str}} T^{\str}
  \]
  is also a contractible operad. From this approach, it would have been harder
  to compute an explicit formula for \(\op^T\), which is needed for example to
  extend the proof assistant \(\catt\) with a meta\-operation computing the
  opposites, and it would not be immediately obvious that the property of being
  free on a computad is preserved by the formation of inverses.
\end{remark}

\subsection{Opposites of hom \texorpdfstring{$\omega$}{ω}-categories}

We will show that the operations of forming hom $\omega$\-categories and
opposite categories commute. To make this statement precise, we first extend the
action of $G = \Z_2^{\N_{>0}}$ on $\omega\Cat$ to an action on bipointed
$\omega$\-categories
\[
  \op : G\to \Aut(\catbip)
\]
by letting for $w\in G$ the opposite of a bipointed $\omega$\-category
$(X,x_,x_+)$ be the $\omega$\-category $\op_w X$ with the same basepoints when
$1\not\in w$, and with the basepoints swapped when $1\in w$.

\begin{lemma}\label{lem:op-susp-comp}
  For every $w\in G$, there exists a natural isomorphism
  \[\op_w^{\susp} : \susp \op_{w-1} \Rightarrow \op_w\susp : \comp\to \compbip\]
  compatible with the natural isomorphism of Lemma~\ref{lem:op-susp} in the
  sense that
  \[\op_w^{\susp}\free = \free \op_w^{\susp}\]
  and the following diagram commutes:
  \[\begin{tikzcd}
    {} & {\susp\op_{w-1}\cell} & {\op_w\susp \cell} \\
    & {\susp\cell \op_{w-1}} & {\op_w\cell\susp} \\
    & {\cell\susp\op_{w-1}} & {\cell\op_w\susp}
    \arrow["{\op_w^{\susp}\cell}", Rightarrow, from=1-2, to=1-3]
    \arrow["{\susp\op_{w-1}^{\cell}}"', Rightarrow, from=1-2, to=2-2]
    \arrow["{\susp^{\cell}\op_{w-1}}"', Rightarrow, from=2-2, to=3-2]
    \arrow["{\op_w\susp^{\cell}}", Rightarrow, from=1-3, to=2-3]
    \arrow["{\op_w^{\cell}\susp}", Rightarrow, from=2-3, to=3-3]
    \arrow["{\cell\op_w^{\susp}}"', Rightarrow, from=3-2, to=3-3]
  \end{tikzcd}\]
\end{lemma}
\begin{proof}
  We will build natural isomorphism
  \[\op_w^{\susp} : \susp \op_{w-1}
    \Rightarrow \op_w\susp : \comp_n\to \comp_{n+1}\]
  inductively on $n \ge -1$ commuting with the forgetful functors $u_{n}$,
  the inclusion functors $\free_n$ and making the following pentagons commute
  \[\begin{tikzcd}
    & {\cell_n} \\
    {\cell_n\op_{w-1}} && {\cell_{n+1}\susp} \\
    {\cell_{n+1}\susp\op_{w-1}} && {\cell_{n+1}\op_w\susp}
    \arrow["{\cell_{n+1}\op_w^{\susp}}", Rightarrow, from=3-1, to=3-3]
    \arrow["{\susp^{\cell}}", Rightarrow, from=1-2, to=2-3]
    \arrow["{\op_w^{\cell}\susp}", Rightarrow, from=2-3, to=3-3]
    \arrow["{\op_{w-1}^{\cell}}"', Rightarrow, from=1-2, to=2-1]
    \arrow["{\susp^{\cell}\op_{w-1}}"', Rightarrow, from=2-1, to=3-1]
  \end{tikzcd}\]
  \[\begin{tikzcd}
    & {\type_n} \\
    {\type_n\op_{w-1}} && {\type_{n+1}\susp} \\
    {\type_{n+1}\susp\op_{w-1}} && {\type_{n+1}\op_w\susp}
    \arrow["{\type_{n+1}\op_w^{\susp}}", Rightarrow, from=3-1, to=3-3]
    \arrow["{\susp^{\type}}", Rightarrow, from=1-2, to=2-3]
    \arrow["{\op_w^{\type}\susp}", Rightarrow, from=2-3, to=3-3]
    \arrow["{\op_{w-1}^{\type}}"', Rightarrow, from=1-2, to=2-1]
    \arrow["{\susp^{\type}\op_{w-1}}"', Rightarrow, from=2-1, to=3-1]
  \end{tikzcd}\]
  For the unique $(-1)$\-computad, we let
  \[\op_{w}^{\susp} : \set{v_-,v_+}\to \set{v_-,v_+}\]
  be the identity function when $1\not\in w$, and the function swapping the two
  generators when $1\in w$. This is a natural isomorphism making the second
  pentagon commute: both sides send the unique $(-1)$\-sphere to the $0$\-sphere
  $(v_-,v_+)$ when $1\not\in w$, and to the $0$\-sphere $(v_+,v_-)$ otherwise.

  For an $n$\-computad $C = (C_{n-1}, V_n^C, \phi_n^C)$ where $n\in\N$, we let
  \[
    \op_{w,C}^{\susp} = (\op_{w,C_{n-1}}^{\susp}, \var) :
      \susp \op_{w-1}C \to \op_w \susp C
  \]
  This is a well-defined morphism of computads by the commutativity of the
  second pentagon one dimension lower. Moreover, it commutes with the forgetful
  and the free functors by construction.

  We will show that the first pentagon commutes for a computad $C$ and that
  $\op_w^{\susp}$ is natural mutually inductively. First we see that the
  pentagon commutes when restricted to generators, since both $\op_w^{\cell}$
  and $\susp^{\cell}$ preserve generators. Suppose now that $n>0$ and let
  $c = \coh{B}{A}{\tau}$ a coherence $n$\-cell of $C$. Then we see that
  \begin{align*}
    (\op_{w,\susp C}^{\cell}\circ& \susp^{\cell}_C)(c) \\
      &= \coh{\op_w\susp B}{A_1}{\op_w\susp\tau\circ\free_n\op_w^{\susp B}}\\
    (\cell_{n+1}(\op_{w,C}^{\susp})\circ
      & \susp^{\type}_{\op_{w-1}C}\circ \op_{w-1,C}^{\type})(c) \\
      &= \coh{\susp\op_{w-1}B}{A_2}
        {\op_{w,C}^{\susp}\circ\susp\op_{w-1}\tau\circ\susp\free_n\op_{w-1}^B}
  \end{align*}
  where
  \begin{align*}
    A_1
      &= \type_n\free_n(\op_w^{\susp B})^{-1}(\op_w^{\type}\susp^{\type}A) \\
    A_2
      &= \susp^{\type}_{\op_{w-1}C}
        (\type_{n-1}\free_{n-1}(\op_{w-1}^B)^{-1}(\op_{w-1}^{\type}A))
  \end{align*}
  By definition of the suspension and the opposite of a tree, we have that
  \[  \op_w\susp B = \susp \op_{w-1}B,\]
  so the trees over which those coherence cells are built agree. Moreover, the
  natural isomorphism $\op_w^{\susp B}$ is defined to be the composite
  \[\op_w^{\susp B} = \op_{w,\Pos{B}}^{\susp} \circ \susp \op_{w-1}^B.\]
  Using that fact and the naturality of \(\susp^{\type}\), we can rewrite the
  spheres $A_1$ and $A_2$ respectively as
  \begin{align*}
    A_1
      &= \type_n\free_n\susp(\op_{w-1}^B)^{-1}
        (\type_n\free_n(\op_{w,\Pos{B}}^{\susp})^{-1}
        (\op_w^{\type}\susp^{\type}A)) \\
    A_2
      &=
        \type_{n}\susp\free_{n-1}(\op_{w-1}^B)^{-1}
        (\susp^{\type}_{\op_{w-1}C}\op_{w-1}^{\type}A)
  \end{align*}
  and observe that they agree by commutativity of the diagram for spheres one
  dimension lower, and commutativity of the suspension with the functor
  $\free_n$. Moreover, we may assume that the naturality square for $\tau$
  commutes by the inductive hypothesis, which shows that the morphisms defining
  the coherence cells agree. Therefore, the first pentagon commutes on coherence
  cells as well.

  Let now $\sigma : D\to C$ be a morphism of $n$\-computads and
  suppose that the pentagon commutes when restricted to cells of the form
  $\sigma_{n,V}(v)$ for $v\in V_n^D$. To show that the naturality square for
  $\sigma$ commutes, i.e. that
  \[
    \op_{w,C}^{\susp} \circ \susp\op_{w-1}\sigma =
    \op_{w}\susp\sigma \circ \op_{w,D}^{\susp},
  \]
  we may assume by induction on the dimension and commutativity with the
  forgetful functors that the underlying morphisms of $n$\-computads agree. It
  remains to show that the two morphisms agree on top-dimensional generators.
  Let therefore $v\in V_{n+1}^{\susp\op_{w-1}D} = V_n^D$ be a generator. Then
  \begin{align*}
    (\op_{w,C}^{\susp} \circ \susp\op_{w-1}\sigma)_V(v)
      &= \cell_{n+1}(\op_{w,C}^{\susp})
        (\susp_{\op_{w-1}C}^{\cell}\op^{\cell}_{w-1,C}(\sigma_V(v))) \\
      &= \op_{w,\susp C}^{\cell}\susp_C^{\cell}(\sigma_V(v)) \\
      &= (\op_w\susp\sigma)_V(v) \\
      &= \cell_{n+1}(\op_w\susp\sigma)(\var v) \\
      &= \cell_{n+1}(\op_w\susp\sigma)((\op_{w,D}^{\susp})_V v) \\
      &= (\op_w\susp\sigma \circ \op_{w,D}^{\susp})_V(v)
  \end{align*}
  so the two morphisms agree on generators as well. Hence, the naturality square
  commutes.

  This concludes the induction on $n\in \N$. By commutativity with the forgetful
  functors, the natural isomorphisms $\op_w^{\susp}$ for every $n\in \N$
  combine to a natural isomorphism
  \[
    \op_w^{\susp} : \susp\op_{w-1} \Rightarrow \op_w\susp : \comp \to \comp
  \]
  as well. Commutativity of the first pentagon shows that the diagram of the
  lemma commutes, since $\op_w^{\susp}\cell$ is the identity on
  positive-dimensional cells.
\end{proof}

\begin{proposition}
  For every $w\in G$, the following diagram commutes
  \[\begin{tikzcd}
    \catbip & \catbip \\
    \omega\Cat & \omega\Cat
    \arrow["\deloop", from=1-2, to=2-2]
    \arrow["\deloop"', from=1-1, to=2-1]
    \arrow["{\op_w}", from=1-1, to=1-2]
    \arrow["{\op_{w-1}}"', from=2-1, to=2-2]
  \end{tikzcd}\]
\end{proposition}

\begin{proof}
  In order to prove commutativity of this diagram for some $w\in G$, it is
  useful to prove commutativity of the analogous diagram on the level of
  globular sets first:
  \[\begin{tikzcd}
    \bipointed & \bipointed \\
    \glob & \glob
    \arrow["\deloop", from=1-2, to=2-2]
    \arrow["\deloop"', from=1-1, to=2-1]
    \arrow["{\op_w}", from=1-1, to=1-2]
    \arrow["{\op_{w-1}}"', from=2-1, to=2-2]
  \end{tikzcd}\]
  The mate of the natural isomorphism of Lemma~\ref{lem:op-susp} is a natural
  transformation fitting in this square defined as the whiskered composite
  \[
    \op_{w-1}\deloop = \deloop\susp\op_{w-1}\deloop
    \xRightarrow{\deloop \op_w^{\susp}\deloop}\deloop \op_w\susp \deloop
    \xRightarrow{\deloop\op_w \kappa} \deloop\op_w
  \]
  for $\kappa$ the counit of the adjunction $\susp\dashv \deloop$. One of the
  snake equations of this adjunction states that $\deloop \kappa$ is an
  identity. Combining that with the fact that $\op_w$ preserves cells and acts
  trivially on morphisms, we see that $\deloop\op_w\kappa$ must also be an
  identity. Moreover, the natural isomorphism $\op_w^{\susp}$ was defined to be
  the identity on positive-dimensional cells, so $\deloop\op_w^{\susp}$ must
  also be an identity. Since the mate of $\op_w^{\susp}$ is an identity natural
  transformation, we conclude that the square above must commute.

  Since the diagram commutes on the level of globular sets, and the forgetful
  functors $U^T$ and $U^{\Tbip}$ are faithful, it follows that the diagram
  commutes on the level of $\omega$\-categories if it commutes for
  objects, meaning that for every bipointed $\omega$\-category $X$, the
  $\omega$\-categories $\deloop\op_w X$ and $\op_{w-1}\deloop X$ are equal.
  Both $\omega$\-categories have the same underlying globular set by
  commutativity of the square on the level of globular sets, so it remains to
  show that they have the same structure morphisms. Unwrapping the definitions
  of $\op_w$ and $\deloop$, this amounts to the commutativity of the following
  diagram of natural transformations
  \[\begin{tikzcd}
    {T\op_{w-1}\deloop} & {\op_{w-1}T\deloop} & {\op_{w-1}\deloop\Tbip} \\
    {T\deloop \op_w} & {\deloop \Tbip \op_w} & {\deloop \op_w\Tbip}
    \arrow["{\op_{w-1}^T\deloop}", Rightarrow, from=1-1, to=1-2]
    \arrow["{\op_{w-1}\deloop^T}", Rightarrow, from=1-2, to=1-3]
    \arrow[equals, from=1-1, to=2-1]
    \arrow[equals, from=1-3, to=2-3]
    \arrow["{\deloop^T\op_w}"', Rightarrow, from=2-1, to=2-2]
    \arrow["{\deloop\op_w^{\Tbip}}"', Rightarrow, from=2-2, to=2-3]
  \end{tikzcd}\]
  where $\op_w^{\Tbip}$ is simply $\op_w^T$ seen as a natural isomorphism
  between bipointed globular sets. By naturality of the mate correspondence,
  commutativity of this diagram is equivalent to that of the following one:
  \[\begin{tikzcd}
    {\susp T \op_{w-1}} & {\susp\op_{w-1}T} & {\op_w\susp T} \\
    {\Tbip\susp\op_{w-1}} & {\Tbip\op_w\susp} & {\op_w\Tbip\susp}
    \arrow["{\susp\op_{w-1}^T}", Rightarrow, from=1-1, to=1-2]
    \arrow["{\op_w^{\susp} T}", Rightarrow, from=1-2, to=1-3]
    \arrow["{\op_w\susp^T}", Rightarrow, from=1-3, to=2-3]
    \arrow["{\susp^T\op_{w-1}}"', Rightarrow, from=1-1, to=2-1]
    \arrow["{\Tbip\op_w^{\susp}}"', Rightarrow, from=2-1, to=2-2]
    \arrow["{\op_w^{\Tbip}\susp}"', Rightarrow, from=2-2, to=2-3]
  \end{tikzcd}\]
  Replacing each $\op_w^T$ by each inverse and rotating the diagram, we are
  left to show that the following diagram commutes:
  \[\begin{tikzcd}[column sep = huge]
    {\susp\op_{w-1}T} & {\susp T\op_{w-1}} & {\Tbip\susp\op_{w-1}} \\
    {\op_w\susp T} & {\op_w\Tbip\susp} & {\Tbip\op_w\susp}
    \arrow["{\susp\op_{w-1}^{\cell}\free}", Rightarrow, from=1-1, to=1-2]
    \arrow["{\susp^{\cell}\free\op_{w-1}}", Rightarrow, from=1-2, to=1-3]
    \arrow["{\Tbip\op_w^{\susp}}", Rightarrow, from=1-3, to=2-3]
    \arrow["{\op_w^{\susp}T}"', Rightarrow, from=1-1, to=2-1]
    \arrow["{\op_w\susp^{\cell}\free}"', Rightarrow, from=2-1, to=2-2]
    \arrow["{\op_w^{\cellbip}\freebip\susp}"', Rightarrow, from=2-2, to=2-3]
  \end{tikzcd}\]
  But this is precisely the diagram of Lemma~\ref{lem:op-susp-comp} whiskered
  on the right with $\free$, hence it commutes.
\end{proof}


%% file: 6-applications.tex
\section{Applications to Eckmann-Hilton cells}
\label{sec:eh}
Throughout this article, we have illustrated the construction of suspensions and
opposites on the various binary composition operations that have introduced in
Section~\ref{sec:cell-compositions}. We will illustrate further our
constructions in action on a more complex construction, that of the
Eckmann-Hilton cells.

\subsection{A clockwise Eckmann-Hilton cell}
The Eckmann-Hilton arguments states that given two monoid structures on a set
satisfying some compatibility relation, the two structures are necessarily equal
and commutative. This arguments can be translated to weak
\(\omega\)\-categories, where it states that considering two cells \(c\) and
\(c'\) whose boundary are both identities over the same cell, there exists a
higher cell witnessing the commutativity of their vertical composite. This cell
comes from the compatibility of the compositions \(\compcell{1}{}{}\) and
\(\compcell{0}{}{}\) on \(1\)\-cells of the \(\omega\)\-category.

We sketch here the definition of such a higher cell, that we call the clockwise
Eckmann-Hilton argument. For this, we consider the computad \(C_{\eh}\) which
has one \(0\)\-dimensional generator \(x\) and two \(2\)\-dimensional generators
\(a,b\) whose sources and targets are all the identity cell \(\id x\) on \(x\).
We will define a \(3\)\-cell \(\eh\) in \(C_{\eh}\) whose source is
\(\compcell{1}{a}{b}\) and whose target it \(\compcell{1}{b}{a}\). Given a weak
\(\omega\)\-category \(\wcat X\), a pair of \(2\)\-cells \(c,c'\) in \(\wcat X\)
whose sources and targets are identities on the same \(0\)\-cell correspond
exactly to a morphism of \(\omega\)\-categories \(f:K^{T}C_{\eh} \to \wcat X\).
The clockwise Eckmann-Hilton cell on \(c\) and \(c'\) in \(\wcat X\) is then
defined to be:
\[
  \eh(c,c') := f(\eh).
\]
The source of this cell is \(\compcell{1}{c}{c'}\) and its target is
\(\compcell{1}{c'}{c}\).

The \(3\)\-cell \(\eh\) can be geometrically illustrated as the following
composite
\[
  \begin{array}{ccccc}
      \begin{tikzcd}[ampersand replacement=\&]
        x
        \ar[r]
        \ar[r, bend right = 25, phantom, "\Downarrow_{b}"]
        \ar[r, bend right = 50]
        \ar[rr, bend left = 80]
        \ar[rr, bend left = 40, phantom, "\Downarrow"]
        \ar[rr, bend right = 80]
        \ar[rr, bend right = 40, phantom, "\Downarrow"]
        \& x
        \ar[r, bend left = 50]
        \ar[r, bend left = 25, phantom, "\Downarrow_{a}"]
        \ar[r]
        \& x
      \end{tikzcd}
    & \Rrightarrow
    &
      \begin{tikzcd}[ampersand replacement=\&]
        x
        \ar[r, bend left]
        \ar[r, phantom, "\Downarrow_{b}"]
        \ar[r, bend right]
        \ar[rr, bend left = 80]
        \ar[rr, bend left = 40, phantom, "\Downarrow"]
        \ar[rr, bend right = 80]
        \ar[rr, bend right = 40, phantom, "\Downarrow"]
        \& x
        \ar[r, bend left]
        \ar[r, phantom, "\Downarrow_{a}"]
        \ar[r, bend right]
        \& x
      \end{tikzcd}
    & \Rrightarrow
    &
      \begin{tikzcd}
        x
        \ar[r]
        \ar[r, bend left = 25, phantom, "\Downarrow_{b}"]
        \ar[r, bend left = 50]
        \ar[rr, bend left = 80]
        \ar[rr, bend left = 40, phantom, "\Downarrow"]
        \ar[rr, bend right = 80]
        \ar[rr, bend right = 40, phantom, "\Downarrow"]
        & x
          \ar[r, bend right = 50]
          \ar[r, bend right = 25, phantom, "\Downarrow_{a}"]
          \ar[r]
        & x
      \end{tikzcd} \\
    \rotatebox{90}{\(\Rrightarrow\)}
    &
    &
    &
    & \rotatebox{-90}{\(\Rrightarrow\)}
    \\
    \begin{tikzcd}[ampersand replacement=\&]
      x
      \ar[rr, bend left = 50]
      \ar[rr, bend left = 25, phantom, "\Downarrow_{a}"]
      \ar[rr]
      \ar[rr, bend right = 25, phantom, "\Downarrow_{b}"]
      \ar[rr, bend right = 50]
      \&\& x
    \end{tikzcd}
    &
    &
    &
    &
      \begin{tikzcd}[ampersand replacement=\&]
        x
        \ar[rr, bend left = 50]
        \ar[rr, bend left = 25, phantom, "\Downarrow_{b}"]
        \ar[rr]
        \ar[rr, bend right = 25, phantom, "\Downarrow_{a}"]
        \ar[rr, bend right = 50]
        \&\& x
      \end{tikzcd}
  \end{array}
\]
where the unlabelled \(1\)\-cells are the identity on \(x\) and the unlabelled
\(2\)\-cells are unitors. The \(3\)\-cells appearing in the diagrams are
composites of unitors, associators and interchangers. The precise definition of
this cell is quite cumbersome and beyond the scope of our paper. It has been
formalised in the proof assistant \catt~\cite{finster_typetheoretical_2017} for
working in finite computads~\cite{benjamin_catt_2024}. We note that the
Eckmann-Hilton cell is not unique, since we can get cells of the same type by
iteratively rotating the cells \(a\) and \(b\) around each other.

\subsection{Opposites of the Eckmann-Hilton cell}
We now discuss the action of the opposite operation on the clockwise
Eckmann-Hilton cell. The first thing that one can notice is that the computad
\(C_{\eh}\) is self-dual, in that \(\op_{w}C_{\eh} = C_{\eh}\) for every
\(w\in G\). Moreover, we have the following equations for the composites of
\(a\) and \(b\):
\begin{align*}
  \op^{\cell}_{1}(\compcell{0}{a}{b}) & = \compcell{0}{b}{a}
  & \op^{\cell}_{1}(\compcell{1}{a}{b}) & = \compcell{1}{a}{b}\\
  \op^{\cell}_{2}(\compcell{0}{a}{b}) & = \compcell{0}{a}{b}
  & \op^{\cell}_{2}(\compcell{1}{a}{b}) & = \compcell{1}{b}{a}.
\end{align*}
It follows that the opposites \(\op_{1}^{\cell}(\eh)\) and
\(\op_{2}^{\cell}(\eh)\) are again Eckmann-Hilton cells, obtained by rotating
\(a\) and \(b\) around each other anticlockwise. They further satisfy the
following equality
\[
  \op_{1}^{\cell}(\eh(a,b)) = \op_{2}^{\cell}(\eh(b,a)),
\]
which we have formally verified with the proof assistant \catt. Given an
\(\omega\)\-category \(\wcat X\) and two cells \(c,c'\) of \(\wcat X\) whose
sources and targets are identities on the same \(0\)\-cell, the anticlockwise
Eckmann-Hilton cell on \(c\) and \(c'\) in \(\wcat X\) can be obtained either as
the clockwise Eckmann-Hilton cell on \(c\) and \(c'\) in the opposite category
\(\op_{1}\wcat X\), or as the clockwise Eckmann-Hilton cell on \(c'\) and \(c\)
in the opposite category \(\op_{2}\wcat X\).

\begin{remark}
  It turns out that the opposite of the Eckmann-Hilton cell \(\op_{1}(\eh)\) is
  also its inverse, as explained in our subsequent
  work~\cite{benjamin_invertible_2024}.
\end{remark}

\subsection{Suspensions of the Eckmann-Hilton cell}
Finally, we discuss the construction of the suspension applied to the
Eckmann-Hilton cell. The computad \(\susp^{n}C_{\eh}\) can be described
explicitly, and it is immediate to see that a morphism of \(\omega\)\-categories
\(K^{T}\susp^{n}C_{\eh} \to \wcat X\) exactly corresponds to a pair of
\((n+2)\)\-cells \(c\) and \(c'\) in \(\wcat X\), whose sources and targets are
all identities over the same \(n\)\-cell. In this situation, we define the
clockwise Eckmann-Hilton cell on \(c\) and \(c'\) to be
\[
  \eh_{n}(c,c') := f((\susp^{\cell})^{n}(\eh)).
\]
This cell witnesses the equivalence between \(\compcell{n+1}{c}{c'}\) and
\(\compcell{n+1}{c'}{c}\), and comes from the interchanger between the
\((n+1)\)\-composition and the \(n\)\-composition of \((n+2)\)\-cells. This
witness cell is obtained as the clockwise Eckmann-Hilton cell on \(c\) and
\(c'\) in the iterated hom \(\omega\)\-category on \(\wcat X\).
